\documentclass[12pt]{article}
\usepackage[latin9]{inputenc}
\usepackage{amsmath,latexsym}
\usepackage[psamsfonts]{amssymb}
\usepackage{url}
\usepackage[bookmarks,bookmarksnumbered,pdfstartview=FitBH]{hyperref}
\hypersetup{
colorlinks=false, 
linkcolor= blue, 
urlcolor= blue,
citecolor=blue,
}

\newtheorem{theoreme}{Th\'eor\`eme}[section]
\newtheorem{proposition}[theoreme]{Proposition}

\newtheorem{definition}[theoreme]{Definition}

\newtheorem{definition-proposition}[theoreme]{Proposition-D\'efinition}

\newtheorem{hypothesis}[theoreme]{Framework}

\newtheorem{theorem}[theoreme]{Theorem}
\newtheorem{lemma}[theoreme]{Lemma}
\newtheorem{corollary}[theoreme]{Corollary}
\newtheorem{example}[theoreme]{Example}

\newtheorem{remark}[theoreme]{Remark}

\newtheorem{question}[theoreme]{Question}

\def\dfn#1{\textbf{\textit{#1}}}  

\def\endofproof{\nobreak\hfill $\blacksquare$ \goodbreak}

\def\Nmath{\mathbb{N}}
\def\Rmath{\mathbb{R}}

\def\Zmath{\mathbb{Z}}

\newcommand{\action}[1]{\ {\curvearrowright}^{#1}}

\newcommand{\ME}{\ {{\overset{\mathrm{ME}}{\sim}}}\ }
\newcommand{\MEa}[1]{\ \underset{\smash{#1}}{{\overset{\mathrm{ME}}{\sim}}}\ } 
 
\def\SOE{\overset{\mathrm{SOE}}{\sim}}
\newcommand{\notME}{\ {{\overset{\mathrm{ME}}{\not\sim}}}\ }

\newcommand{\locus}[1]{{\mathrm{D}}(#1)}

\def\OE{\overset{\mathrm{OE}}{\sim}}

\def\freeprod{\mathop{*}}

\def\Image{\text{\rm Im }}

\def\FF{{\bf F}}

\def\RR{{\mathcal{R}}}
\def\SS{{\mathcal S}}
\def\VV{\mathcal{V}}
\def\TT{\mathcal{T}}

\newcommand{\gauche}[1]{{#1}}
\newcommand{\droite}[1]{{#1}'}

\def\kk{k}
\def\KK{K}
\def\qq{q}
\def\EE{{\droite{\mathcal{E}}}}
\def\WW{\mathcal{W}}
\def\UUU{\droite{U}}
\def\VVV{\droite{V}}
\def\WWW{\droite{W}}

\def\gi{\gauche{p}}
\def\gI{\gauche{P}}
\def\di{\droite{p}}
\def\dI{\droite{P}}

\def\FI{$\mathcal{FI}$}
\def\MFI{$\mathcal{MFI}$}
\def\SST{\widetilde{\Sigma}}

\begin{document}

\thispagestyle{empty}
\title{Free products, Orbit Equivalence and Measure Equivalence Rigidity}
\author{Aur\'elien Alvarez and Damien Gaboriau\thanks{CNRS} }
\date{\today}
\maketitle

\begin{abstract}
{We study the analogue in orbit equivalence of free product decomposition and free indecomposability for countable groups. We introduce the (orbit equivalence invariant) notion of freely indecomposable ({\FI}) standard probability measure preserving equivalence relations and establish a criterion to check it, namely non-hyperfiniteness and vanishing of the first $L^2$-Betti number. We obtain Bass-Serre rigidity results, \textit{i.e.} forms of uniqueness in free product decompositions of equivalence relations with ({\FI}) components. 
The main features of our work are weak algebraic assumptions and no ergodicity hypothesis for the components.
We deduce, for instance, that a measure equivalence between two free products of non-amenable groups with vanishing first $\ell^2$-Betti numbers is induced by measure equivalences of the components. We also deduce new classification results in Orbit Equivalence and II$_1$ factors.
}
\end{abstract}

\section{Introduction}

Bass-Serre theory \cite{Ser77} studies groups acting on trees and offers extremely powerful tools to understand their structure, together with a geometric point of view that illuminates several classical results on free product decompositions. For instance Kurosh's subgroup theorem \cite{Kurosh-1934}, that describes the subgroups in a free product of groups and, as a by-product, the essential uniqueness in free product decompositions into freely indecomposable subgroups, is much easier to handle via Bass-Serre theory.

In Orbit Equivalence theory, the notion of free products or freely independent standard equivalence relations introduced in \cite{Gab00a} proved to be useful in studying the cost of equivalence relations and for some classification problems. The purpose of our article is connected with the uniqueness condition in free product decompositions, in the measurable context. To this end, we will take full advantage of the recent work of the first named author \cite{Alv08a, Alv08b}, who develops a Bass-Serre theory in this context. In particular, Theorem~\ref{th: a la Kurosh restrict} and Theorem~\ref{th: a la Kurosh} will be crucial for our purpose.

\medskip

Very roughly, the kind of results we are after claim that if a standard measured equivalence relation is decomposed in two ways into a free product of factors that are not further decomposable in an appropriate sense, then the factors are pairwise related.
However, due to a great flexibility in decomposability, it appears that certain types of free decomposition, namely slidings (Definition~\ref{def: sliding}) and slicings (Definition~\ref{def: slicing}),  are banal and somehow inessential (see Section~\ref{subsect: slicings-slidings}).
We thus start by clearing up the notion of a \dfn{freely indecomposable ({\FI})} standard measured countable equivalence relation (Definition~\ref{def: freely indec}), ruling out inessential decompositions (Definition~\ref{def: trivialization of decompositon}).

A countable group $\Gamma$ is said \dfn{measurably freely indecomposable ({\MFI})} if all its free probability measure preserving (p.m.p.) actions produce freely indecomposable ({\FI}) equivalence relations.
As expected, a free product of two infinite groups is not {\MFI}, and in fact none of its free p.m.p. actions is {\FI}. The same holds for infinite amenable groups (\textit{cf.} Corollary~\ref{Cor: act infinite amenable non-FI}). 
On the other hand, freely indecomposable groups in the classical sense are not necessarily {\MFI}, for instance the fundamental group of a compact surface of genus $\geq 2$ (see Proposition~\ref{prop: MFI is ME invariant}). 
We now give a prototypical instance of our results:
\begin{theorem}\label{th: MFI groups, ergodic restrict}
Consider two families of infinite countable {\MFI} groups $(\Gamma_i)_{i\in I}$ and $(\Lambda_j)_{j\in J}$, $I=\{1,2,\cdots, n\}$, $J=\{1,2,\cdots, m\}$, $n,m\in \Nmath^{*}\cup\{\infty\}$. Consider two free probability measure preserving actions $\alpha$ and $\beta$ of the free products
  on standard Borel spaces whose restrictions to the factors $\alpha\vert \Gamma_i$ and $\beta\vert \Lambda_j$ are ergodic.
If the actions $\alpha$ and $\beta$ are stably orbit equivalent  
\begin{equation}
(\freeprod\limits_{i\in I}\Gamma_i)\action{\alpha}(X,\mu)\  \SOE\ (\freeprod\limits_{j\in J}\Lambda_j)\action{\beta} (Y,\nu)
\end{equation}
 then $n=m$ and there is a bijection $\theta:I\to J$ for which the restrictions are stably orbit equivalent
\begin{equation}
\alpha\vert{\Gamma_i} \SOE \beta\vert{\Lambda_{\theta(i)}}
\end{equation}
\end{theorem}
Of course, such a statement urges us to exhibit {\MFI} groups,
and it appears that their class is quite large:
\begin{theorem}[Cor. \ref{cor: beta-1=0 -> Gamma is MFI}]
Every non-amenable countable group $\Gamma$ with vanishing first $\ell^2$-Betti number ($\beta_1(\Gamma)=0$) is measurably freely indecomposable ({\MFI}).
\end{theorem}
Recall that the $\ell^2$-Betti numbers are a sequence of numbers $\beta_r(\Gamma)$ defined by Cheeger-Gromov \cite{CG86} attached to every countable discrete group $\Gamma$ and that they have a general tendency to concentrate in a single dimension $r$ and to vanish in the other ones (see \cite{BV97}, \cite{Luc02}). The first $\ell^2$-Betti number vanishes for many "usual" groups, for instance for amenable groups, direct products of infinite groups, lattices in $\textrm{SO}(p,q)$ ($p.q\not=2$), lattices in $\textrm{SU}(p,q)$, groups with Kazhdan's property (T). It is worth noting that infinite Kazhdan's property (T) groups also follow {\MFI}  from Adams-Spatzier \cite[Th 1.1]{AS90} (see \cite[Ex. IV.12]{Gab00a}).
The list of groups with vanishing $\beta_1$ may be continued, for instance, with the groups with an infinite finitely generated normal subgroup of infinite index, groups with an infinite normal subgroup with the relative property (T),  amalgamated free products of groups with $\beta_1=0$ over an infinite subgroup, mapping class groups, ...
On the other hand, for a free product of two (non trivial) groups we have $\beta_1(\Gamma_1*\Gamma_2)>0$ unless $\Gamma_1=\Gamma_2=\Zmath/2\Zmath$, in which case $\Gamma_1*\Gamma_2$ is amenable.

\medskip
Results in the spirit of Theorem~\ref{th: MFI groups, ergodic restrict} were obtained as by-products of operator algebraic considerations in \cite[Cor. 0.5, Cor. 7.6, Cor. 7.6']{IPP05}, and also very recently in \cite[Cor. 6.7]{CH08}.  Our results cover a large part of these corollaries.
We will come back more precisely on the differences between these papers and ours, but an important issue is that they both require the ergodicity of the actions restricted to the factors and some particular algebraic assumptions on the groups.

\medskip
We will extend our framework by introducing marginal free groups or relatives (recall that in Kurosh's theorem there are "vertex subgroups" and a free group) and more seriously by removing the ergodicity assumption on the actions of the factors; and both of these extensions prove to be necessary to handle with Measure Equivalence of groups (see \cite{Gab05a} for a survey on this notion introduced by M. Gromov). 
Recall that two countable groups $\Gamma$ and $\Lambda$ are \dfn{Measure Equivalent (ME)}, in symbols:
\begin{equation}
\Gamma\MEa{\kappa}\Lambda
\end{equation}
if and only if they admit Stably Orbit Equivalent (SOE) free p.m.p. actions.
The real number $\kappa\in \Rmath^{*}_{+}$ is called the \dfn{generalized index} or the \dfn{compression constant} according to whether one focuses on the classification of groups up to ME or on more operator algebraic aspects of Orbit Equivalence. Commensurable groups are ME, and the generalized index then coincides with the usual index for subgroups. 
It is proved in~\cite[{\textbf{P}$_{\!\mbox{\tiny{ME}}}\mathbf{6}$}, p. 1814-1816]{Gab05a}) that measure equivalent groups with generalized index $1$ induce measure equivalence of their free products: 
\begin{center}
if $\Gamma_i\MEa{1}\Lambda_i$ then $\freeprod\limits_{i\in I}\Gamma_i\MEa{1}\freeprod\limits_{i\in I}\Lambda_i$.
\end{center}
 Our technics allow us to settle a converse when the factors are {\MFI}.
We observe that being {\MFI} is a Measure Equivalence invariant (Proposition~\ref{prop: MFI is ME invariant}), and before stating our ME result, we consider the following striking example. It prevents us from being overoptimistic or expecting a bijective correspondence between the factors.
\begin{example}\label{ex: prevent overoptimism}
If $\Gamma'_1\triangleleft \Gamma_1$ and $\Gamma'_2 \triangleleft\Gamma_2$ are two normal subgroups of finite index $\kappa$ such that $\Gamma_1/\Gamma'_1\simeq \Gamma_2/\Gamma'_2\simeq K$, then the following groups are ME with $\Gamma_1*\Gamma_2$ with generalized index $\kappa$ and they satisfy:
\begin{equation}
\Gamma'_1*\underbrace{\Gamma_2*\Gamma_2*\cdots*\Gamma_2}_{\kappa \mathrm{\ copies}}\MEa{1}\underbrace{\Gamma_1*\Gamma_1*\cdots*\Gamma_1}_{\kappa \mathrm{\ copies}}*\Gamma'_2\MEa{1}\Gamma'_1*\Gamma'_2*\FF_{\kappa-1}
\end{equation}
where $\FF_p$ is the free group on $p$ generators.
In fact these three groups are even mutually commensurable with finite kernels and generalized index $1$, since they appear as 
the kernels of the three natural epimorphisms $\Gamma_1*\Gamma_2 \twoheadrightarrow K$.
\end{example}
Theorem~\ref{th: MFI groups, ergodic restrict} ensures that such "pathologies" are ruled out by adding ergodic assumptions on the actions of the factors.
Explicit actions witnessing these measure equivalences are easily produced by suspension, and the fact that they are not ergodic when restricted to some factors is not at all incidental (and the above Example~\ref{ex: prevent overoptimism} may be better understood). We are able, from Theorem~\ref{th: rigidity BS}, to localize some constraints on the failure of ergodicity, for instance for any action witnessing a measure equivalence between the following commensurable groups:
\begin{corollary}\label{cor: SOE -> non erg on a restrict} Assume that $\Gamma_1, \Gamma_2$ are {\MFI} and not ME, and that $\Gamma'_1$ has finite index $\kappa\geq 2$ in $\Gamma_1$.
Then, for any stably orbit equivalent actions $\Gamma_1*\Gamma_2\action{\alpha} X$ and $\Gamma'_1*{\Gamma_2*\Gamma_2*\cdots*\Gamma_2}\action{\beta}Y$, the restriction $\alpha\vert \Gamma_2$ is not ergodic.
\end{corollary}
We are now in position to state our general Measure Equivalence result:
\begin{theorem}[ME Bass-Serre rigidity]\label{th: ME BS rigidity-intro}
Consider two families of infinite countable
\- {\MFI} groups 
	(for instance non-amenable with vanishing first $\ell^2$-Betti number) 
	$(\Gamma_i)_{i\in I}$ and $(\Lambda_j)_{j\in J}$ 
	, $I=\{1,2,\cdots, n\}$, $J=\{1,2,\cdots, m\}$, $n,m\in \Nmath^{*}\cup\{\infty\}$. If their free products are measure equivalent, 
\begin{equation}
\freeprod\limits_{i\in I}\Gamma_i\ME \freeprod\limits_{j\in J}\Lambda_j
\end{equation}
 then there are two maps $\theta : I\to J$ and $\droite{\theta}:J\to I$
such that:
\begin{eqnarray}
\Gamma_i \ME \Lambda_{\theta(i)} \textrm{\ \ and \ \ }
\Lambda_j \ME  \Gamma_{\droite{\theta}(j)}
\end{eqnarray}
Moreover, if $\Gamma_0,\Lambda_0$ are two groups in the ME classes of some free groups, then the same conclusion holds under the assumption:
 \begin{equation}\label{eq: ME between free prod}
\freeprod\limits_{i\in I}\Gamma_i\freeprod \Gamma_0\ME \freeprod\limits_{j\in J}\Lambda_j\freeprod \Lambda_0.
\end{equation}
\end{theorem}
Observe that we do not assume the generalized index $\kappa=1$. Would we do so, we would not get $\kappa=1$ in the conclusion as Example~\ref{ex: prevent overoptimism} again indicates.
Also observe that the groups $\Gamma_0, \Lambda_0$ do not appear in the conclusion.


It is interesting to observe that one may combine our theorem~\ref{th: ME BS rigidity-intro} for free products with Monod-Shalom's Theorem~1.16 \cite{MS06} for direct products. Using the facts that a free product of infinite groups belongs to their class $C_{reg}$, and that a direct product of non-amenable  groups is {\MFI}, we get the following type of results: 
\begin{corollary}
Assume $\Gamma$ is either \\
i) a finite direct product of non-trivial free products of torsion-free {\MFI} groups $\Gamma_i$, or \\
ii) a free product of non-trivial finite direct products of torsion free groups $\Gamma_i$ in the class $C_{reg}$.
\\
If $\Gamma$ is measure equivalent with a group $\Lambda$ of the same kind, then the elementary pieces $\Gamma_i$ of $\Gamma$ define the same set of ME-classes as those of $\Lambda$. 
\end{corollary}
Of course, this contruction can be iterated by taking alternatively free or finite direct products of groups $\Gamma_i$ as in i) or ii) above, according to whether the first operation is a free or a direct product.
A measure equivalence with a group $\Lambda$ of the same kind entails measure equivalences between the elementary pieces. Notice that the number of iterations follows the same for $\Gamma$ and $\Lambda$.

The above Theorem~\ref{th: ME BS rigidity-intro} is essentially a consequence of the following SOE Theorem~\ref{th: SOE BS-rigidity} (see Theorem~\ref{th: rigidity BS} for the more general measured equivalence relations statement).
We continue with similar data. 
In the SOE context, the role of free groups is played by treeability (see \cite{Ada88,Gab00a} or Section~\ref{subsect: graphings-treeings} for more on this notion).
Recall that it follows from \cite{Hjo06} that a group $\Gamma$ is ME with a free group if and only if it admits a free p.m.p. treeable action (see \cite[\textbf{P}$_{\!\mbox{\tiny{ME}}}\mathbf{8}$]{Gab05a}). A group is said \dfn{strongly treeable} if all its free p.m.p. actions are treeable. This is for instance the case of the amenable groups (even finite or even trivial) or the free products of amenable groups (ex. free groups). We have no example of a group that 
is  ME with a free group but that is not strongly treeable.

We display separately the following framework that will be in use in the next results:
\begin{hypothesis}\label{hyp: general framework}
Let\\
-- $(\gauche{\Gamma}_{\gi})_{\gi\in \gI}$ and $(\droite{\Gamma}_{\di})_{\di\in \dI}$ 
	, $\gI=\{1,2,\cdots, \gauche{n}\}$, 	 $\dI=\{1,2,\cdots, \droite{n}\}$, $\gauche{n},\droite{n}\in \Nmath^{*}\cup\{\infty\}$ be two families of infi\-nite countable groups
\\
-- $\gauche{\Gamma}_0, \droite{\Gamma}_0$ be two countable groups.\\
	Assume $\gauche{\alpha}$ and $\droite{\alpha}$ are two free p.m.p. stably orbit equivalent actions 
	\begin{equation}
	(\freeprod\limits_{\gi\in \gI}\gauche{\Gamma}_{\gi}\freeprod \gauche{\Gamma}_0)
	\action{\gauche{\alpha}}(\gauche{X},\gauche{\mu})\ \ \SOE\ \ (\freeprod\limits_{\di\in \dI}
	\droite{\Gamma}_{\di}\freeprod \droite{\Gamma}_0)\action{\droite{\alpha}}(\droite{X},
	\droite{\mu})
\end{equation}
of the free products on standard Borel spaces,  such that \\
	 -- the restrictions $\gauche{\alpha}\vert \gauche{\Gamma}_{\gi}$ and $\droite{\alpha}\vert \droite{\Gamma}_{\di}$ are freely indecomposable ({\FI}) (for $\gi,\di\in \gI,\dI$),\\
-- the restrictions $\gauche{\alpha}\vert {\gauche\Gamma}_0$ and $\droite{\alpha}\vert \droite{\Gamma}_0$ are treeable.
\end{hypothesis}

\begin{theorem}[SOE Bass-Serre Rigidity]\label{th: SOE BS-rigidity} 
If $\gauche{\alpha}$ and $\droite{\alpha}$ are two SOE actions
 as in Framework~\ref{hyp: general framework}, then 
up to countable partitions, the components are in one-to-one correspondence in the following sense. There exist
\begin{enumerate}
	\item 
	\label{item: left partition} for each $\gi\in \gI$, a measurable $\gauche{\alpha}\vert\gauche{\Gamma}_{\gi}$-invariant partition $\gauche{X}=\coprod\limits_{\gauche{k}\in \gauche{K}(\gi)} \gauche{X}_{\gauche{k}}$
	\item 
	\label{item: right partition} 
	for each $\di\in \dI$, a measurable $\droite{\alpha}\vert\droite{\Gamma}_{\di}$-invariant partition $\droite{X}=\coprod\limits_{\droite{k}\in \droite{K}(\di)} \droite{X}_{\droite{k}}$	
	\item a bijection $\theta: \coprod\limits_{{\gi}\in {\gI}}\gauche{\KK}(\gi)\to \coprod\limits_{{\di}\in {\dI}}\droite{\KK}(\di)$ between the index sets
\end{enumerate}
according to which the restrictions of the actions to the factors and the subsets are SOE:
\begin{equation}
\forall\gauche{k}\in \coprod\limits_{{\gi}\in {\gI}}\gauche{\KK}(\gi)
\hskip 40pt
\gauche{\alpha}\vert \gauche{\Gamma}_{\gauche{k}}\times \gauche{X}_{\gauche{k}}\SOE \droite{\alpha}\vert \droite{\Gamma}_{\theta(\gauche{k})}\times \droite{X}_{\theta(\gauche{k})} \hskip 40pt
\end{equation}
with the obvious notational conventions: $\gauche{\Gamma}_{\gauche{k}}:=\gauche{\Gamma}_{\gi}$ for the unique $\gi\in \gI$ such that $\gauche{k}\in \gauche{K}(\gi)$, and $\droite{\Gamma}_{\theta(\gauche{k})}:=\droite{\Gamma}_{\di}$ for the unique $\di\in \dI$ such that $\theta(\gauche{k})\in \droite{K}(\di)$.
\end{theorem}
Observe that under ergodic assumptions on the actions restricted to the factors, the invariant partitions turn trivial and $\theta$ gives a bijective correspondence between the original index sets. Ergodicity on one side may force the same situation:
\begin{corollary}\label{Cor: ergodicity on one side} 
Consider two SOE actions $\gauche{\alpha}$ and $\droite{\alpha}$ as in Framework~\ref{hyp: general framework}. Assume that $\gauche{n}\leq \droite{n} <\infty$ and that the restrictions of the actions to the $\gauche{\Gamma}_{\gi}$-factors  $\gauche{\alpha}\vert\gauche{\Gamma}_{\gi}$ are ergodic, $\forall \gi\in \gI$.
Then the restrictions $\droite{\alpha}\vert \droite{\Gamma}_{\di}$ are also ergodic, $\gauche{n}= \droite{n}$ and $\theta$ gives a bijection $\theta:\gI\to \dI$.
If moreover $\beta_1(\gauche{\Gamma}_{\gi})=\beta_1(\droite{\Gamma}_{\di})=0$ for all $\gi,\di \in \gI,\dI$ and $\gauche{\Gamma}_0=\droite{\Gamma}_0=\{1\}$, then the factors follow orbit equivalent 
$\gauche{\alpha}\vert\gauche{\Gamma}_{\gi}\OE \droite{\alpha}\vert \droite{\Gamma}_{\theta(\di)}$.
\end{corollary}

We also get some consequences for ergodic components from Theorem~\ref{th: rigidity BS}.
It is a banal observation that the number of ergodic components (let's denote it $\# \mathrm{erg\ comp}
(\sigma)$) of a single action  $G\action{\sigma}(X,\mu)$ is invariant under stable orbit equivalence. We obtain a survival of this invariant for a restriction of an action to factors of a free product:
\begin{corollary} Consider two SOE actions $\gauche{\alpha}$ and $\droite{\alpha}$ as in Framework~\ref{hyp: general framework}. Let $\gI_{1}\subset \gI$ and $\dI_{1}\subset \dI$ be the indices of those groups $\gauche{\Gamma}_{\gi}, \droite{\Gamma}_{\di}$ that are measure equivalent with $\gauche{\Gamma}_1$. Then the number of ergodic components of the restrictions are equal:
\begin{equation}
\sum_{\gi\in \gI_{1}} \# \mathrm{erg\ comp}(\gauche{\alpha}\vert \gauche{\Gamma}_{\gi})=
\sum_{\di\in \dI_{1}} \# \mathrm{erg\ comp}(\droite{\alpha}\vert \droite{\Gamma}_{\di})
\end{equation}
\end{corollary}
Even, the measures of the ergodic components become a SOE invariant, under a control of the self generalized indices. 
The set $I_{\mathrm{ME}}(\Gamma)$ of possible generalized indices~$\kappa$ in measure equivalences between a group $\Gamma$ and itself $\Gamma\MEa{\kappa} \Gamma$ is an invariant of the ME class of $\Gamma$ (see \cite{Gab02b} or \cite[\textbf{P}$_{\!\mbox{\tiny{ME}}}\mathbf{17}$]{Gab05a}).
The condition $I_{\mathrm{ME}}(\Gamma)=\{1\}$ is obtained for instance when $\Gamma$ has an $\ell^2$-Betti number $\beta_q(\Gamma)\not=0, \infty$ \cite{Gab02}.
For sake of simplicity, we give a sample of the kind of statements that may be derived from Theorem~\ref{th: rigidity BS} (see Theorem~\ref{th: IME trivial -> OE}):
\begin{corollary}\label{cor: beta-p>0 --> OE}
Assume that the $(\gauche{\Gamma}_{\gi})_{\gi\in \gI}$ have vanishing $\beta_1$ and that $\gauche{\Gamma}_0$ is a free group. Assume that $\gauche{\Gamma}_1$ admits at least one $\ell^{2}$-Betti number $\beta_q(\gauche{\Gamma}_1)$ different from $0$ and $\infty$ and that it is not measure equivalent with any of the other $\gauche{\Gamma}_{\gi}$, $\gi\not=1$. If $\Theta$ is a SOE between
two  p.m.p. actions
$\gauche{\alpha}$ and $\droite{\alpha}$ of $(\freeprod_{\gi\in \gI}\gauche{\Gamma}_{\gi}\freeprod \gauche{\Gamma}_0)$, then $\Theta$ is in fact an OE and the restrictions to $\gauche{\Gamma}_1$ are OE. In particular, they have the same measure space of ergodic components.
\end{corollary}
\begin{corollary}
Let $\Gamma_0=\FF_2$ and $\Gamma_1=\FF_3\times \FF_3$. Consider a one-parameter family of free p.m.p. action $\Gamma_0\freeprod \Gamma_1\action{\alpha_s}(X,\mu)$, where the restriction $\alpha_s\vert \Gamma_1$
has two ergodic components of respective measures $s, 1-s$. The actions $\alpha_s$ are not mutually stably orbit equivalent for $s\in [0,1/2]$.
\end{corollary}

\medskip

Recall that free p.m.p. group actions $\Gamma\action{\sigma}(X,\mu)$ define finite von Neumann algebras by the so called group-measure space construction of Murray-von Neumann or von Neumann crossed product $L^{\infty}(X,\mu)\rtimes_{\sigma} \Gamma$. 
Stably orbit equivalent actions define stably isomorphic crossed-products, but the converse does not hold in general, and this leads to the following definition. Two free p.m.p. actions $\gauche{\Gamma}\action{\gauche{\sigma}}(\gauche{X}, \gauche{\mu})$ and  $\droite{\Gamma}\action{\droite{\sigma}}(\droite{X}, \droite{\mu})$ are called \dfn{von Neumann stably equivalent} if there is $\kappa\in (0,\infty)$ such that $L^{\infty}(\gauche{X}, \gauche{\mu})\rtimes_{\gauche{\sigma}} \gauche{\Gamma}\simeq (L^{\infty}(\droite{X}, \droite{\mu})\rtimes_{\droite{\sigma}} \droite{\Gamma})^{\kappa}$.

\medskip

Both papers \cite{IPP05}, \cite{CH08} establish rigidity phenomena in operator algebras and derive orbit equivalence results for the components of free products from an assumption of von Neumann stable equivalence on the actions. To this end, some strong algebraic constraints on the involved groups are imposed. More precisely in \cite[Cor. 0.5, Cor. 7.6, Cor. 7.6']{IPP05}, the analysis relies on the notion of relative property (T) in von Neumann algebras introduced by S.~Popa in \cite{Pop06}, and thus 
the groups $\gauche{\Gamma}_{\gi}, \droite{\Gamma}_{\gi}$ (in the notation of Framework~\ref{hyp: general framework}) are required to admit a non virtually abelian subgroup with the relative property (T) and some ICC-like and normal-like properties (for instance, they may be ICC property (T) groups) (see \cite[Assumption 7.5.1]{IPP05}). In \cite[Cor. 6.7]{CH08}, the operator algebraic notion involved is primality, so that the assumption on the groups $\gauche{\Gamma}_{\gi}, \droite{\Gamma}_{\gi}$ is to be ICC non-amenable direct products of infinite groups. 
In both cases, they all satisfy $\beta_1=0$.
As already mentioned, the actions restricted to the factors $\gauche{\alpha}\vert\gauche{\Gamma}_{\gi}$ and $\droite{\alpha}\vert \droite{\Gamma}_{\di}$ are assumed to be ergodic. On the other hand, the assumption that the actions are SOE $\gauche{\alpha}\SOE\droite{\alpha}$ is replaced by the weaker one that $\gauche{\alpha}$ and $\droite{\alpha}$ are von Neumann stably equivalent.
All these results exploit the antagonism between free products and either various forms of property (T), or direct product, or more generally in our case the vanishing of the first $\ell^2$-Betti number. 
Also, in \cite{IPP05} the "marginal" groups $\gauche{\Gamma}_0, \droite{\Gamma}_0$ are solely assumed to be a-T-menable (\textit{i.e.} to have Haagerup property) a property that in turn is antagonist to property (T), while our $\gauche{\Gamma}_0, \droite{\Gamma}_0$ are ME with a free group (antagonist to {\FI}) and thus a-T-menable. Observe that there are {\MFI} a-T-menable groups, thus able to play the role of a $\gauche{\Gamma}_{\gi}$ or a $\gauche{\Gamma}_0$ according to the approach.

\bigskip

On the other hand, it follows from  \cite[Th. 7.12, Cor.7.13]{IPP05} that von Neumann stable equivalence entails stable orbit equivalence, among the free p.m.p. actions 
of free products of (at least two) infinite groups, as soon as one of the two actions has the relative property (T) in the sense of \cite[Def. 4.1]{Pop06}.
Meanwhile, Theorem 1.2 of \cite{Gab08} establishes that any free product of at least two infinite groups admits a continuum of relative property (T) von Neumann stably inequivalent ergodic free p.m.p. actions, whose restriction to each free product component is conjugate with any prescribed (possibly non-ergodic) action.

When injected in our context, this gives further classifications results for II$_1$ factors. For instance:
\begin{theorem}
Let $\gauche{\Gamma}_{1}, \gauche{\Gamma}_{2}$ be non-ME, non-amenable groups with $\beta_1=0$.
Assume $\beta_q(\gauche{\Gamma}_{1})\not=0,\infty$ for some $q>1$. 
The crossed-product II$_1$ factors $M_1\freeprod_{A} M_2$ associated with the various ergodic relative property (T) free p.m.p. actions $\gauche{\Gamma}_{1}\freeprod \gauche{\Gamma}_{2}\action{\sigma}(X,\mu)$ are classified by the pairs $A\subset M_1$, and in particular by the isomorphism class of the centers $Z(M_1)$ of the crossed-product associated with the restriction of the action to $\gauche{\Gamma}_1$,
equipped with the induced trace. 
\end{theorem}
Of course, we do not claim that this invariant is complete.

\bigskip

Our treatment considers p.m.p. standard equivalence relations instead of just free p.m.p. group actions. The notion of $L^2$-Betti numbers for these objects, introduced in \cite{Gab02}, gives a criterion for free indecomposability:
\begin{theorem}[Th.~\ref{th: rel beta-1=0 +nowhere amen ->FI}]
If $\RR$ is a nowhere hyperfinite p.m.p. standard equivalence relation on $(X,\mu)$ with $\beta_{1}(\RR, \mu)=0$, then it is freely indecomposable.
\end{theorem}

Our main result is Theorem~\ref{th: rigidity BS} which describes the kind of uniqueness one can expect in a free product decomposition into {\FI} subrelations.

\bigskip

Some parts of our work may be led in the purely Borel theoretic context.
For instance, we show that a treeable {\FI} equivalence relation is necessarily smooth (Proposition~\ref{rem: treeable + FI -> smooth}). 
Theorem~\ref{th: ME BS rigidity-intro}, Corollary~\ref{Cor: ergodicity on one side} and Corollary~\ref{cor: beta-p>0 --> OE} are proved in Section~\ref{sect: proof of the Corollaries}.

\section{Free Product Decompositions}
\subsection{Generalities}

Let $X$ be a standard Borel space. All the equivalence relations we will consider are Borel with countable classes. By countable, we mean "at most countable". In the measured context, $X$ is equipped with a non-atomic finite measure $\mu$ and the equivalence relations are measure preserving (m.p.) (\textit{resp.} probability measure preserving (p.m.p.) when $\mu$ is a probability measure) and the following definitions are understood up to a null-set.

Since we are about to consider, on a standard Borel space $X$, equivalence relations that may be defined only on a subset of $X$, 
we set:
\begin{definition}\index{domain}
The Borel set on which a countable standard Borel equivalence relation $\RR$ is defined will be called its 
\dfn{domain} and will be denoted by $\locus{\RR}$.
\end{definition}
Recall that a \dfn{complete section} is a Borel subset of $\locus{\RR}$ that meets all the classes.
A \dfn{fundamental domain} is a Borel subset of $\locus{\RR}$ that meets each class  exactly once. 
An equivalence relation is \dfn{smooth} if it admits a fundamental domain. An equivalence relation is \dfn{finite} if all its classes are finite. In this case, it is smooth. In the probability measure preserving case, $\RR$ is smooth if and only if it is finite.
$\RR$ is \dfn{aperiodic} if its classes are all infinite on $\locus{\RR}$.
If $U\subset \locus{\RR}$, we denote by $\RR\vert U$ the \dfn{restriction} $\RR\cap U\times U$ of $\RR$ to $U$, and its domain is $\locus{\RR\vert U}=U$.
The relation $\RR$ is \dfn{trivial} if $\RR=\{(x,x): x\in \locus{\RR}\}$, \textit{i.e.} if its classes are reduced to singletons. The equivalence relation $\RR$ on $\locus{\RR}\subset X$ naturally \dfn{extends} to an equivalence relation on the whole of $X$ by setting the classes of $x\in X\setminus \locus{\RR}$ to be reduced to the singleton $\{x\}$. We use the same notation $\RR$ for the extended relation since it will be clear from the context what we are considering.
The \dfn{full group} $[\RR]$ of $\RR$ is the group of all Borel isomorphisms of $\locus{\RR}$ whose graph is contained in $\RR$. The \dfn{full pseudogroup} $[[\RR]]$ is the family of all Borel partial isomorphisms between Borel subsets of $\locus{\RR}$ whose graph is contained in $\RR$. The equivalence relation  $\SS$ is a \dfn{subrelation} of $\RR$ if $\locus{\SS}\subset \locus{\RR}$ and $(x,y)\in \SS$ implies $(x,y)\in \RR$. If $\SS$ is a subrelation of $\RR$ and $\phi:A\to B $ is a partial isomorphism in the full pseudogroup $[[\RR]]$ of $\RR$ whose target $B$ is contained in $\locus{\SS}$ then 
\begin{equation}
\phi^{-1} \SS \phi
\end{equation}
 denotes the equivalence relation of domain $A$ defined by $(x,y)\in \phi^{-1} \SS \phi$ if and only if $(\phi(x),\phi(y))\in \SS$. It is the image of $\SS\vert B$ under $\phi^{-1}$. Two subrelations $\SS_1$ and $\SS_2$ of $\RR$ are said \dfn{inner conjugate} in $\RR$ if there is a partial isomorphism $\phi\in [[\RR]]$ with domain $\locus{\SS_2}$ and target $\locus{\SS_1}$ such that $\SS_2= \phi^{-1} \SS_1 \phi$.

\subsection{Free Products}

\begin{definition}[{see \cite[D\'ef.~IV.9]{Gab00a}}]
A countable family of equivalence relations $(\RR_i)_{i\in I}$ with domains $\locus{\RR_i}\subset X$ is \dfn{freely independent} if the following holds:
for any $n$-tuple $(x_1, \cdots, x_n)$ of elements of $X$ such that $x_n=x_1$ and $(x_j,x_{j+1})\in \RR_{i_j}$ there is an index $j$ such that $i_j=i_{j+1}$.
The equivalence relation $\RR$ is \dfn{decomposed as the free product}
\begin{equation}
\RR=\freeprod_{i\in I} \RR_i
\end{equation}
or \dfn{is the free product} of the countable family $(\RR_i)_{i\in I}$
if the family of subrelations  is freely independent and generates $\RR$ (in particular $\locus{\RR}=\cup_{i\in I} \locus{\RR_i}$). The $\RR_i$ are the \dfn{factors} or the \dfn{components} of the free product decomposition.
\end{definition}

\begin{lemma}\label{lem: intersection de conj de facteurs}
Let $\RR=\RR_1*\RR_2*\RR_3$ be decomposed as a free product. Consider $\SS_1$ and $\SS_2$ two subrelations of $\RR_1$ and $\RR_2$ that are inner conjugate in $\RR$, then $\SS_1$ (and $\SS_2$) are smooth.
\end{lemma}
Proof: Let $\phi\in [[\RR]]$, $\phi:\locus{\SS_2}\to \locus{\SS_1}$ such that $\SS_2=\phi^{-1} \SS_1 \phi$. Assume first that $\phi$ decomposes as a product of partial isomorphisms taken \textit{strictly} from the $[[\RR_i]]$, \textit{i.e.} $\phi= \phi_{r_n}\cdots \phi_{r_2}\phi_{r_1}$, such that for each $j$: $\phi_{r_j}\in [[\RR_{k_j}]]$, $k_{j}\not=k_{j+1}$ and for every $z$ in its domain $\phi_{r_j}(z)\not=z$. Any $(x,y)\in \SS_2$ defines, by introducing the right subwords of $\phi$, a "rectangular" cycle
\begin{eqnarray*}
x & \overset{\RR_{k_1}}{\sim} \phi_{r_1}(x)\overset{\RR_{k_2}}{\sim} \phi_{r_2}\phi_{r_1}(x)\sim\cdots \overset{\RR_{k_{n-1}}}{\sim}\phi_{r_{n-1}}\cdots \phi_{r_2}\phi_{r_1}(x) \overset{\RR_{k_n}}{\sim}&\phi(x)\\
{}_{\SS_2}{\wr} & & \ \ {\wr}_{\SS_1}\\
y &\overset{\RR_{k_1}}{\sim} \phi_{r_1}(y)\overset{\RR_{k_2}}{\sim} \phi_{r_2}\phi_{r_1}(y)\sim\cdots \overset{\RR_{k_{n-1}}}{\sim}\phi_{r_{n-1}}\cdots \phi_{r_2}\phi_{r_1}(y)\overset{\RR_{k_n}}{\sim}&\phi(y)
\end{eqnarray*}
that may be shorten by definition of free products. Due to strictness, the only possible  shortenings may occur around the vertical sides: after a possible shortening of the horizontal sides in case $k_1=2$ or $k_n=1$,  the extreme points have to coincide, so that $\SS_1$ and $\SS_2$ are trivial. 
The general case reduces to this after a decomposition of the domain of $\phi$ into pieces where it satisfies the above assumption. Its restrictions to the pieces being trivial, $\SS_1$ and $\SS_2$ follow smooth.
\endofproof

\subsection{Graphings and Treeings}\label{subsect: graphings-treeings}

Recall from \cite{Lev95}, \cite{Gab00a} that a countable family of partial isomorphisms $\Phi=(\phi_i)_{i\in I}$ is called a \dfn{graphing} and defines an equivalence relation $\RR_{\Phi}=\langle \Phi\rangle = \langle \phi_i: i\in I\rangle $ on $\locus{\RR_{\Phi}}=$ the union of the domains and the targets of the $\phi_i$'s.
It is a \dfn{treeing} if any equation $\phi_{i_1}^{\varepsilon_1}\phi_{i_2}^{\varepsilon_2}\cdots \phi_{i_n}^{\varepsilon_n}(x)=x$ (with $\varepsilon_{i_j}=\pm 1$) implies there is an index $j$ such that $i_j=i_{j+1}$ and $\varepsilon_{i_j}=-\varepsilon_{i_{j+1}}$.
An equivalence relation is \dfn{treeable} if it admits a generating treeing.
The notion of treeing were introduced by S.~Adams \cite{Ada88} and proved to be very useful in \cite{Gab98} and \cite{Gab00a}.

We recall some properties of treeable equivalence relations and their connections with free products.
\begin{proposition}\label{prop: properties of treeable equiv. rel.} The following holds:
\begin{enumerate}
\item A subrelation of a treeable equivalence relation is itself treeable. \label{item: subrel treeable}
\item If $\Phi=(\phi_i)_{i\in I}$ is a treeing, then $\RR_{\Phi}$ is the free product of the subrelations generated by the individual partial isomorphisms $\RR_{\Phi}=\freeprod_{i\in I} \langle \phi_i \rangle$.\label{item: treeing and free product}
\item A free product of treeable equivalence relations is treeable.\label{item: free prod of treeable is treeable}
\item A treeable equivalence relation is freely decomposed as a free product of finite subrelations. \label{item: treeable : free prod finite subrel.}

\end{enumerate}
\end{proposition}
Proof:
Item \ref{item: subrel treeable} is Theorem IV.4 in \cite{Gab00a} (where the proof does not make use of the measure). This has also been shown independently by Jackson-Kechris-Louveau \cite{JKL02}. See also \cite{Alv08a} for a geometric approach. Item \ref{item: treeing and free product} is immediate from the definitions (see \cite[Ex. IV.10]{Gab00a}). So is also item \ref{item: free prod of treeable is treeable}: a treeing for the free product is made of the union of treeings for the factors. 
As for item~\ref{item: treeable : free prod finite subrel.}, 
by a result of Slaman-Steel and Weiss (\cite{SS88}, \cite{Wei84}), each singly generated relation $\langle \phi_i \rangle$ is hyperfinite. 
We now claim that a hyperfinite equivalence relation is a free product of finite equivalence relations. To that end, it is enough to show that if ${\cal R}_1$ is a subrelation of a finite equivalence relation ${\cal R}_2$ then there exists a (finite) subrelation ${\cal R}'_1$ of ${\cal R}_2$ such that ${\cal R}_2={\cal R}_1 \freeprod {\cal R}'_1$. Indeed, given fundamental domains $D_1$ and $D_2$ of ${\cal R}_1$ and ${\cal R}_2$, we get a Borel finite-to-one projection $\pi:D_1 \longrightarrow D_2$ whose pre-images naturally define the required ${\cal R}'_1$ on $D_1$.
\endofproof

\medskip
An \dfn{action} $\Gamma\action{}(X,\mu)$ is \dfn{treeable} if the equivalence relation it generates is treeable. If $\Gamma$ is measure equivalent with a free group then it admits a treeable p.m.p. free action. 
\begin{question}
Are there groups with both treeable and non-treeable free p.m.p. actions~?
\end{question}

\subsection{Slidings and Slicings}\label{subsect: slicings-slidings}

We now consider two banal ways of freely decomposing an equivalence relation.
\begin{definition}[Slicing]\label{def: slicing}\index{slicing}
The \dfn{slicing}  \dfn{affiliated} with a  Borel partition $\coprod_{j\in J} V_j$ of the domain $\locus{\RR}$ into $\RR$-invariant subsets,
is the free product decomposition 
\begin{equation}
\RR=\freeprod\limits_{j\in J} \RR\vert V_j.
\end{equation}
\end{definition}
\begin{definition}[Sliding]\label{def: sliding}\index{sliding}
Let $U\subset\locus{\RR}$ be a complete section of $\RR$. A \dfn{sliding} of $\RR$ to $U$ consists in a smooth treeable subrelation $\TT<\RR$ defined on $\locus{\RR}$ with fundamental domain $U$ and in the corresponding free product decomposition 
\begin{equation}
\RR=\RR{\vert U}*\TT.
\end{equation}
\end{definition}
An explicit construction of such a treeable subrelation $\TT$ for each such $U$ may be found in \cite[Lem. II.8]{Gab00a}, where the notion is introduced in a measured context. Notice that the proof does not make use of the measure. 
From this, one can deduce an easy particular case of Theorem~\ref{th: a la Kurosh restrict} from next section:
\begin{proposition}
If $\RR= \freeprod_{j\in J} \RR_j$ is a free product decomposition with $\locus{\RR_{j_0}}=\locus{\RR}$ and $U$ is a complete section for each $\RR_j$, then $\RR{\vert U}= \freeprod_{j\in J} \RR_{j}{\vert U} \freeprod \TT$, where $\TT$ is a treeable subrelation.
\end{proposition}
Proof: Consider slidings $\RR_j=\RR_j\vert U\freeprod \TT_j$ and inject them in the decomposition of $\RR= \freeprod_{j\in J} (\RR_{j}{\vert U}\freeprod\TT_j)$. 
For $j\not= {j_0}$, the sliding $\TT_j\freeprod \TT_{j_0}=(\TT_j\freeprod \TT_{j_0})\vert U \freeprod \TT_{j_0}$, gives the global sliding  $$\freeprod\limits_{j\in J}\TT_j=\freeprod\limits_{j\in J\setminus\{{j_0}\}} (\TT_j\freeprod \TT_{j_0})\vert U\freeprod \TT_{j_0}.$$ It follows that 
$$\RR\vert U= \freeprod_{j\in J} \RR_{j}{\vert U} \freeprod \overbrace{\freeprod_{j\in J\setminus\{{j_0}\}} (\TT_j\freeprod \TT_{j_0})\vert U}^{\TT:=}$$
where $\TT$ is treeable by Proposition~\ref{prop: properties of treeable equiv. rel.}, items~\ref{item: subrel treeable} and \ref{item: free prod of treeable is treeable}. \endofproof

\section{\texorpdfstring{Theorems à la Kurosh after \cite{Alv08a}}{Theorems à la Kurosh after Alv08}}

We will make a crucial use in our construction of some tools introduced by the first named author, namely the following two analogues of Kurosh's theorem \cite{Kurosh-1934} for subgroups of free products, in the context of p.m.p. standard equivalence relations. The first one concerns the particular situation of a subrelation which is simply the restriction to some Borel subset of a given free product. 
\begin{theorem}[A la Kurosh for restrictions \cite{Alv08a}]\label{th: a la Kurosh restrict}
Let 
\begin{equation}
	\SS=\freeprod\limits_{i\in I} \SS_i
\end{equation}
be a free product decomposition of $\SS$ and $Y\subset X$ a complete section for $\SS$. Then $\SS$ admits a refined free product decomposition induced by slicings of the factors $\SS_i$
\begin{equation}
	\SS_i=\freeprod\limits_{k\in K(i)} \SS_i\vert X_k \hskip 40pt D(\SS_i)=\coprod_{k\in K(i)} X_k 
\end{equation}
such that the restriction $\SS{\vert Y}$ admits a free product decomposition:
\begin{equation}
	\SS{\vert Y}=\freeprod\limits_{i\in I}\bigl(\freeprod\limits_{k\in K(i)} \VV_k  \bigr)*\TT
\end{equation}
where $\TT$ is a treeable subrelation; and for each $i\in I$:
\begin{enumerate}

\item for each $k\in K(i)$,
 there is a partial isomorphism $\phi_k\in [[\SS]]$, defined on the domain $\locus{\VV_k}$,
 that inner conjugates $\VV_k$ with  $\SS_{i}$ restricted to the target of $\phi_k$: 
\begin{equation}
\VV_k=\phi_{k}^{-1}\ \SS_{i}\  \phi_k
\end{equation}
and $X_k=\SS_i \phi_k(\locus{\VV_k})$ is the $\SS_i$-saturation of the image $\phi_k(\locus{\VV_k})$;

\item \label{item: S-i cap Y is part of Kurosh decomposition}
if $\locus{\SS_i}\cap Y$ is non-negligible, then there is $k\in K(i)$ such that  $\VV_k=\SS_{i}{\vert \locus{\SS_i}\cap Y}$ (\textit{i.e.} $\locus{\VV_k}=\locus{\SS_i}\cap Y$, $\phi_k=id_{\locus{\SS_i}\cap Y}$).
\label{it: exists S-i for any i}
\end{enumerate}

\end {theorem}

Compare with the analog result \cite[Prop. 7.4 2$^\circ$]{IPP05}; both the statement and the proof are much less intricate, due to the assumption that the factors $\SS_i$ are ergodic.
This Theorem~\ref{th: a la Kurosh restrict} is itself of course a little bit more precise than the next one which describes the general situation of a subrelation in a free product.

\begin{theorem}[A la Kurosh \cite{Alv08a}]\label{th: a la Kurosh}
Let 
\begin{equation}
	\SS=\freeprod\limits_{i\in I} \SS_i
\end{equation} 
be a free product decomposition of $\SS$.
If $\RR<\SS$ is a subrelation of $\SS$ with non-null domain  $\locus{\RR}\subset X$, then $\RR$ admits a free product decomposition
\begin{equation}
	\RR=\freeprod\limits_{i\in I}\bigl(\freeprod\limits_{k\in K(i)} \VV_k  \bigr)*\TT
\end{equation}
where $\TT$ is a treeable subrelation; and for each $i\in I$:
\begin{enumerate} 
\item for each $k\in K(i)$,
 there is a partial isomorphism $\phi_k\in [[\SS]]$ defined on the domain $\locus{\VV_k}$ such that 
\begin{equation}
	\VV_k=\RR\cap \phi_{k}^{-1}\ \SS_{i}\  \phi_k
\end{equation}
In particular, $\VV_k$ is inner conjugate with a subrelation of $\SS_{i}$.
\item \label{item: R cap S-i is part of Kurosh decomposition}
there is $k\in K(i)$ such that $\VV_k=\RR\cap  \SS_{i}$ (when this intersection is not trivial), $\locus{\VV_k}=\locus{\RR}\cap \locus{\SS_{i}}$ and $\phi_k=id_{\locus{\VV_k}}$.
\end{enumerate}
\end {theorem}

\begin{remark} 
It does not matter whether $\locus{\RR}$ is a complete section of $\SS$ or not.
\end{remark}

\begin{remark}\label{rem: group the treeable part}
If one of the original factors, say $\SS_{i_0}$, is treeable, then the factors $\VV_k$ associated with $k\in K(i_0)$ in Theorem~\ref{th: a la Kurosh restrict}
 or Theorem~\ref{th: a la Kurosh} follow treeable  so that $\freeprod_{k\in K(i_0)} \VV_k\freeprod \TT$ is treeable (see Proposition~\ref{prop: properties of treeable equiv. rel.}). More generally, if the $\VV_k$ in a certain collection 
are treeable, one may assemble them together with $\TT$ to form a treeable relation that may be put in place of $\TT$, in the above theorems.
\end{remark}

\section{Freely indecomposable equivalence relations}

Observing that slicings (Definition~\ref{def: slicing}) and slidings (Definition~\ref{def: sliding})  decompose an equivalence relation as a free product in a somewhat trivial way leads to set the following definition.
\subsection{Free Indecomposability}
\begin{definition}[Inessential Free Product Decomposition]
\label{def: trivialization of decompositon}\index{inessential free product decomposition}
A free product decomposition $\RR=*_{j\in J} \RR_j$, of a (countable) standard Borel equivalence relation on $X$ 
is called \dfn{inessential} if there is a Borel set $U\subset X$ such that:
	\begin{enumerate}
		\item $U$ admits a Borel $\RR{\vert U}$-invariant partition $U=\coprod_{j\in J} U_j$;
		\item for each $j\in J$,  $\RR{\vert U_j}=\RR_{j}{\vert U_j}$; \label{def: triv decomp. pt restriction}
		\item $U$ is a complete section for $\RR$.
	\end{enumerate}
We then say that the partition is a \dfn{trivialization} of the decomposition.
\end{definition}
\begin{remark}\label{Rem: identities up to null sets} In the measured context, all the identities are understood up to a set of measure zero.
\end{remark}

\begin{remark} \label{Rem: sur la def de decomp. triviale}
If the decomposition $\RR=\freeprod_{j\in J} \RR_j$ is trivialized by $U=\coprod_{j\in J} U_j$:
	
	\begin{enumerate}
	\item 
It induces the slicing $\RR{\vert U}=\freeprod_{j\in J} \RR_{j}{\vert U_j}= \freeprod_{j\in J} \RR{\vert U_j}$, and $U$ being a complete section, a sliding/slicing decomposition $\RR=\freeprod_{j\in J} \RR_{j}{\vert U_j}\freeprod\TT$,
 where $\TT$ is a smooth treeable equivalence relation with fundamental domain $U$. 
	\item \label{item: slicing induced by trivialization}
	We insist that the $\RR{\vert U}$-invariance of the partition means that the $\RR$-saturations $V_j:=\RR U_j$ of the $U_j$ have trivial mutual intersections and partition $X$, leading to the slicing:
	\begin{equation}\label{eq: slicing induced by trivialization}
		\RR=\freeprod_{j\in J} \RR\vert V_j
	\end{equation}
\end{enumerate}
\end{remark}
\begin{proposition}\label{Prop: sur la def de decomp. triviale}
Assume the free product decomposition $\RR=\freeprod_{j\in J} \RR_j$ is trivialized by $U=\coprod_{j\in J} U_j$. Then,
\begin{enumerate}
	\item \label{item: R-j satur U-j trivialize}
	$\overline{U}=\coprod_{j\in J} \overline{U}_j$, where $\overline{U}_j:=\RR_j U_j$ is the $\RR_j$-saturation of $U_j$, also trivializes the free product decomposition.
	
	\item 	\label{item: R-j restricted complement bar(U-i) is smooth}
	If $j\in I$\\
	-- the subrelation $\RR_j$ is trivial when restricted to 
	$\overline{U_i}$, 
	for $i\not=j$.\\
	-- the subrelation $\RR_j$ is smooth when restricted to $X\setminus \overline{U_j}$.
	
	\item \label{item: R-i nowhere-smooth-> R=R-i} If $\RR$ is measure preserving, $\locus{\RR_i}=\locus{\RR}$ and $\RR_i$ is aperiodic, then $X=V_i$  and $\RR=\RR_i$ almost everywhere.
	
	\item If $\RR$ is ergodic, $U$ equals one of the $U_j$'s, say $U_{j_0}$, and $\RR=\RR_{j_0}{\vert U}*\TT$, where $\TT$ is a smooth treeing admitting $U$ as fundamental domain. 
\end{enumerate}
\end{proposition}
Proof: The only (maybe) non-obvious facts are items~\ref{item: R-j satur U-j trivialize} and \ref{item: R-j restricted complement bar(U-i) is smooth}:
\\
\ref{item: R-j satur U-j trivialize}.
If $\overline{x}, \overline{y} \in \overline{U}_j$ are $\RR$-equivalent, there are $\RR$-equivalent points $x,y\in U_j$ such that $x\RR_j\overline{x}$ and $ \overline{y} \RR_j y$. By $\RR\vert U_j=\RR_j\vert U_j$, $\overline{x}, \overline{y}$ follow $\RR_j$-equivalent.
\\
\ref{item: R-j restricted complement bar(U-i) is smooth}.
The first part is clear. As for the second part, 
	it is then enough to show (since $X=\cup_j V_j$) that $\RR_j$ restricted to $V_i\setminus \overline{U_i}=\RR U_i \setminus \RR_i U_i$ is smooth for every $i\in I$.
The $\RR_i$-slicing affiliated with $V_i=\overline{U_i}\coprod V_i\setminus \overline{U_i}$ and the $\RR$-invariance of $V_i$ lead to the free product decomposition
$\RR\vert V_i=\RR_i\vert \overline{U_i}\freeprod \RR_i\vert (V_i\setminus \overline{U_i})\freeprod(\freeprod_{j\in J\setminus\{i\}} \RR_j\vert V_i)$. Any partial isomorphism $\phi\in [[\RR]]$ with domain $\subset \overline{U_i}$ and target $\subset V_i\setminus \overline{U_i}$ inner conjugates the restriction of each of the other factors with a subrelation of $\RR\vert \overline{U_i}=\RR_i\vert \overline{U_i}$. By Lemma~\ref{lem: intersection de conj de facteurs}, these subrelations are smooth. Since $\overline{U_i}$ is a complete section for $\RR\vert V_i$, there are enough such $\phi$, and the conclusion follows.\endofproof
\begin{definition}[{\FI} Equivalence Relation]
\label{def: freely indec} \index{freely indecomposable}\index{{\FI}}
A countable standard Borel equivalence relation $\RR$ is \dfn{freely indecomposable} ({\FI}) if any free product decomposition $\RR=\RR_{1}*\RR_{2}*\cdots*\RR_{i}*\cdots$ is inessential in the sense of Definition~\ref{def: trivialization of decompositon}.
\end{definition}
For instance, any finite equivalence relation is {\FI}. 
See Remark~\ref{Rem: identities up to null sets} in the measured context.
\begin{proposition}\label{rem: treeable + FI -> smooth}
If $\RR$ is treeable (for instance hyperfinite) and {\FI} then it is smooth.
In fact, any p.m.p. aperiodic treeable equivalence relation admits an essential decomposition in {\em two} pieces $\RR=\RR_1\ast \RR_2$.
\end{proposition}
\begin{proposition}\label{prop: def FI with 2 pieces}
For p.m.p. standard aperiodic equivalence relations, the Definition~\ref{def: freely indec} can be stated equivalently for free product decompositions in {\em two} pieces.
\end{proposition}
Proof of Prop.~\ref{rem: treeable + FI -> smooth}:
From Proposition~\ref{prop: properties of treeable equiv. rel.}, $\RR$ decomposes as a free product $\freeprod_{i\in I}\RR_i$ of finite subrelations. The property  $\RR{\vert U_i}=\RR_{i}{\vert U_i}$ in Definition~\ref{def: trivialization of decompositon}(\ref{def: triv decomp. pt restriction}) implies that $\RR$ is smooth. 
In the p.m.p. context, the $X$ splits in two $\RR$-invariant Borel subsets $X_{\infty}\coprod X_h$ where the classes have infinitely many ends (resp. where $\RR$ is hyperfinite) \cite{Ada90}.
On $X_h$, $\RR$ is generated by a free action of $\Zmath/2\Zmath\ast \Zmath/2\Zmath$. On $X_{\infty}$, one can find a treeing $\Phi=(\varphi_1)\vee \Phi_2$ of $\RR$ such that the domain $A_1$ of $\varphi_1$ is non-negligible and $\Phi_2$ generates a subrelation $\RR_2$ on $X_{\infty}$ that is aperiodic. The free product decomposition $\RR\vert X_{\infty}=\langle\varphi_1\rangle \ast \RR_2 $ is essential: Prop.~\ref{Prop: sur la def de decomp. triviale}(\ref{item: R-j restricted complement bar(U-i) is smooth}) with aperiodicity of $\RR_2$ implies $\overline{U_2}=X$ (a.e.) for any hypothetical trivialization. But again by Prop.~\ref{Prop: sur la def de decomp. triviale}(\ref{item: R-j restricted complement bar(U-i) is smooth}), $\langle\varphi_1\rangle$ should be trivial on $\overline{U_2}=X$.
\endofproof

\medskip

Proof of Prop.~\ref{prop: def FI with 2 pieces}:
Let $\mathcal{R}$ be a standard aperiodic p.m.p. equivalence relation.
Assuming that any free product decomposition into two pieces $\mathcal{R}=\mathcal{S}_1 \ast \mathcal{S}_2$ is inessential, we show that any free product decomposition $\mathcal{R} = \ast_{j \in J}\mathcal{R}_j$ into at most countably many pieces is also inessential.

  For each $i\in J$, let $\mathcal{R}'_i= \ast_{j \in J \setminus \{i\}}  \mathcal{R}_j$ and consider a trivialization $U_i \coprod U'_i$ associated with the free product decomposition, $\mathcal{R}=\mathcal{R}_i \ast \mathcal{R}'_i$. Moreover $U_i$ may be assumed $\RR_i$-saturated (Prop.~\ref{Prop: sur la def de decomp. triviale} (\ref{item: R-j satur U-j trivialize})).
We claim that $V=\coprod U_i$ trivializes the free product decomposition $\ast_{j \in J}  \mathcal{R}_j$. Since $\RR_i\vert U_i=\RR\vert U_i$ is aperiodic (when $U_i$ is non negligible) and $\RR_i\vert X\setminus U_i$ is smooth (Prop.~\ref{Prop: sur la def de decomp. triviale} (\ref{item: R-j restricted complement bar(U-i) is smooth})), the only point to check is that $V$ is a fundamental domain. The complement of its saturation $Y=X \setminus \mathcal{R}{V}$ is contained in $\cap_{i\in J} (X\setminus U_i)$, so that $\RR\vert Y=\ast_{j \in J}\mathcal{R}_j\vert Y$ is treeable. 

Under a p.m.p. asumption,  Proposition~\ref{rem: treeable + FI -> smooth} and $\mu(Y)\not=0$ would produce a essential free product decomposition in two pieces of $\RR\vert Y$ (thus also of $\RR$), leading to a contradiction. Thus $Y$ is negligible.
\endofproof

\medskip
By Ornstein-Weiss theorem \cite{OW80}, Prop.~\ref{rem: treeable + FI -> smooth} entails:
\begin{corollary}\label{Cor: act infinite amenable non-FI}
Every free p.m.p. action of an infinite amenable group is non-{\FI}.
\end{corollary}
\subsection{Properties of Free Indecomposability}
\begin{remark}\label{rem: extension remains FI}
Given an $\RR$-invariant partition $\locus{\RR}=Y\coprod Z$, then $\RR$ is {\FI} if and only if $\RR\vert Y$ and $\RR\vert Z$ are {\FI}. In particular, if $\RR$ is {\FI}, then the extension of $\RR$ to $X$ (by trivial classes outside $\locus{\RR}$) is also {\FI}. 
\end{remark}

\begin{proposition}\label{prop: restr of FI is FI}
If $\RR$ is {\FI}, then for every non-null Borel set $Y$, the restriction $\RR{\vert Y}$ is also {\FI}.
\end{proposition}
Proof: By Remark~\ref{rem: extension remains FI}, one may assume that $Y$ is a complete section of $\RR$ and that $\RR$ is aperiodic. Any free product decomposition $\RR{\vert Y} = \freeprod_{i \in I} \RR_i$ leads by sliding to $\RR = \freeprod_{i \in I} \RR_i \freeprod \TT$ where $\TT$ is a treeing with fundamental domain $Y$. The {\FI}-property for $\RR$ gives a trivialization $\coprod_{i\in I} U_i \coprod U_{\TT}$, such that $U_i \subset \locus{\RR_i}\subset Y$, $\RR{\vert U_i} = \RR_i{\vert U_i}$ and  $\RR{\vert U_{\TT}} = \TT{\vert U_{\TT}}$. 
Since $\TT$ is smooth, $U_{\TT}$ is negligible (by aperiodicity). It follows that $\coprod_{i\in I} U_i\subset Y$ gives a trivialization for the restriction of $\RR\vert Y$.
\endofproof
\begin{proposition}[Stable Orbit Equivalence Invariance]\label{prop: Stable Orbit Equivalence Invariance}
If $\RR$ and $\SS$ are stably orbit equivalent, then $\RR$ is {\FI} if and only if $\SS$ is {\FI}.
\end{proposition}
The statement being clear for Orbit Equivalence, it remains to show:
If $Y$ a complete section of $\RR$ and $\RR{\vert Y}$ is {\FI}, then $\RR$ is also {\FI}. One may assume that $\RR$ is aperiodic. For each free product decomposition $\RR = \freeprod_{i \in I} \RR_i$, Theorem~\ref{th: a la Kurosh restrict} delivers a free product decomposition  of the restriction $\RR{\vert Y}=\freeprod_{i\in I}\bigl(\freeprod_{k\in K(i)} \phi_{k}^{-1}\ \RR_{i}\  \phi_k  \bigr)*\TT$, with $\TT$ treeable.
The {\FI}-property for $\RR\vert Y$ gives a trivialization $\coprod_{i\in I}\bigl(\coprod_{k\in K(i)} U_k) \coprod U_{\TT}$, where in particular
$(\RR{\vert Y}){\vert U_{\TT}} = \TT{\vert U_{\TT}}$ is treeable (Proposition~\ref{prop: properties of treeable equiv. rel.}~item~\ref{item: subrel treeable}) and thus smooth (by Propositions~\ref{rem: treeable + FI -> smooth} and~\ref{prop: restr of FI is FI}). It follows that $U_{\TT}$ is negligible (by aperiodicity), and from $(\RR\vert Y)\vert U_k=(\phi_{k}^{-1}\ \RR_{i}\  \phi_k)\vert U_k=\phi_{k}^{-1}\ \RR_{i}\vert \phi_k(U_k)\  \phi_k$ that $\coprod_{i\in I}\bigl(\coprod_{k\in K(i)} \phi_k(U_k))$ trivializes the original decomposition.\endofproof
\begin{definition}
A countable group is called \dfn{measurably freely indecomposable} ({\MFI}) if all its free p.m.p. actions are freely indecomposable ({\FI}).
\end{definition}
\begin{proposition}\label{prop: MFI is ME invariant} Being measurably freely indecomposable is a measure equivalence invariant: if $\Gamma\ME \Lambda$ then $\Gamma$ is {\MFI} iff $\Lambda$ is {\MFI}.
\end{proposition}
For instance, such groups as the fundamental group of a compact surface of genus $\geq 2$ are freely indecomposable in the classical sense (they have only one end) but are not {\MFI} since being ME with a non-cyclic free group.

\medskip
\noindent
{Proof of Proposition~\ref{prop: MFI is ME invariant}:}
\def\RRT{\widetilde{\RR}}
Consider two standard  p.m.p. equivalence relations $\RRT$ on $(\tilde{X},\tilde{\mu})$ and $\RR$ on $(X,\mu)$. Let $p:\tilde{X}\to X$ be a measurable map
such that $p_*(\tilde{\mu})\sim\mu$
and  $p$ induces, for (almost) every $\tilde{x}\in \tilde{X}$, a bijection between the $\RRT$-class of $\tilde{x}$ and the $\RR$-class of $p(\tilde{x})$. Such a $p$ is a \dfn{locally bijective morphism from $\RRT$ to $\RR$}.
\begin{remark}
This notion were introduced in \cite[p.1815]{Gab05a} as locally one-to-one and onto morphism from $\RRT$ to $\RR$, and we take this opportunity to correct a regrettable translation mistake that led to use the words one-to-one instead of bijective, all along the paper.
\end{remark}

\begin{lemma}
If $\RRT$ is {\FI} then $\RR$ is also {\FI}.
\end{lemma}
Proof of the lemma: Observe that for a Borel subset $\tilde{W}\subset\tilde{X}$, $\RRT\vert \tilde{W}$ is smooth if and only if $\RR\vert p(\tilde{W})$ is smooth (the restrictions are smooth iff their saturations are smooth, iff their classes are finite).
Let $\RR=*_{j\in J} \RR_j$ be a free product decomposition of $\RR$.
It induces via $p$ a free product decomposition accordingly $\RRT=*_{j\in J} \RRT_j$, where $p$ becomes a locally bijective morphism from $\RRT_j$ to $\RR_j$ (see \cite{Gab05a}). Let $\coprod_{j\in J} \tilde{U}_j$ be a trivializing partition of this decomposition such that the $\tilde{U}_j$ are $\RRT_j$ saturated (Proposition~\ref{Prop: sur la def de decomp. triviale} item \ref{item: R-j satur U-j trivialize}) and let $U_j=p(\tilde{U}_j)$. If $\RRT_j$ is smooth when restricted to some Borel subset $\tilde{Z}_j\subset\tilde{U}_j$, then $\RRT$ follows smooth on $\RRT \tilde{Z}_j$ (by $\RRT_j\vert \tilde{U}_j=\RRT\vert \tilde{U}_j$) and the same holds for $\RR$ on $p(\RRT \tilde{Z}_j)$, so that $\RR$ is {\FI} on this saturated part. One thus may assume that the $\RRT_j\vert \tilde{U}_j$-classes are all infinite. By Proposition~\ref{Prop: sur la def de decomp. triviale} item~\ref{item: R-j restricted complement bar(U-i) is smooth}, $\RRT_j$ is smooth outside $\tilde{U}_j$, so that $p(\tilde{X}\setminus \tilde{U}_j)\cap p(\tilde{U}_j)$ is negligible and $\tilde{U}_j=p^{-1}p(\tilde{U}_j)$. It follows that $\coprod p(\tilde{U}_j)$ is a trivializing partition for $\RR=*_{j\in J} \RR_j$.\endofproof

From this lemma, one gets that if $\Gamma\ME \Lambda$ and if $\Gamma$ is not {\MFI} (\textit{i.e.} $\Gamma$ admits some non {\FI} p.m.p. free action $\Gamma\action{\alpha}(X,\mu)$) then there is a p.m.p. free action of $\Gamma$ that is both non-{\FI} and SOE with a p.m.p. free action of $\Lambda$. 
Let $(\Omega,\nu)$ be a measure equivalence coupling between $\Gamma$ and $\Lambda$.
Consider the coupling $(\Omega\times X,\nu\times \mu)$, with the diagonal actions induced from $\Gamma\action{\alpha}(X,\mu)$ and the trivial action of $\Lambda$ on $X$. The quotient actions $\Gamma\action{} (\Omega\times X)/\Lambda$
and $\Lambda\action{} \Gamma\backslash(\Omega\times X)$ are free (see \cite{Gab02b}), SOE and the first one is non-{\FI} by the above lemma, since it factors onto $\Gamma\action{\alpha}(X,\mu)$. The conclusion of Proposition~\ref{prop: MFI is ME invariant} then follows by SOE invariance (Proposition~\ref{prop: Stable Orbit Equivalence Invariance}).\endofproof

\begin{question}
Are there groups that admit some {\FI} and some non-{\FI} free p.m.p. actions ?
\end{question}

\subsection{\texorpdfstring{$L^2$-Betti numbers}{L2-Betti numbers}}

We now consider finite-measure preserving equivalence relations.
\begin{definition}
A measure preserving standard equivalence relation $\RR$  on $(X,\mu)$ is called \dfn{nowhere hyperfinite} if for every non-null Borel subset $V\subset X$, the restriction $\RR{\vert V}$ is not hyperfinite.
\end{definition}
We establish a criterion for equivalence relations to be {\FI}. The notion of $L^2$-Betti numbers is introduced in \cite{Gab02}. Some useful properties are recalled in Section~\ref{subsect: preliminaries}.
\begin{theorem}\label{th: rel beta-1=0 +nowhere amen ->FI}
If $\RR$ be a nowhere hyperfinite finite-measure preserving standard equivalence relation on $(X,\mu)$ with $\beta_{1}(\RR)=0$, then it is freely indecomposable.
\end{theorem}

In case $\locus{\RR}\not=X$ and $\beta_1(\RR)=0$ (of course computed with respect to the restriction of the measure to $\locus{\RR}$ -- see Section~\ref{subsect: preliminaries}), then the extension of $\RR$ to $X$ (by trivial classes outside $\locus{\RR}$) is also {\FI} (see Remark~\ref{rem: extension remains FI}).

We will prove (Section~\ref{sect: proof th: trivial decomposition for beta-1=0}) more precisely:
\begin{theorem}\label{th: trivial decomposition for beta-1=0}
Assume $\RR$ is an aperiodic finite measure preserving standard equivalence relation with $\beta_{1}(\RR)=0$ and assume that $\RR$ decomposes as a free product $\RR=\RR_{1}*\RR_{2}*\cdots*\RR_{i}*\cdots$. Let $U_{i}\subset X$ be the union of the infinite $\RR_{i}$-classes, for $i=1,2,\cdots$. Then \\
(1) the mutual intersections are trivial, $\mu(U_{i}\cap U_{j})=0$ for $i\not=j$;\\
(2) the partition $U=U_1\coprod U_2\coprod\cdots\coprod U_i\coprod\cdots$ is $\RR{\vert U}$-invariant;\\
(3) the restrictions $\RR_{i}{\vert U_j}$ are trivial for $i\not=j$;\\
(4) $\RR{\vert U_i}=\RR_{i}{\vert U_i}$ for each $i$;\\
(5) the restriction of $\RR$ to the complement $W$ of the saturation $\RR U$ is hyperfinite (if $\mu(W)>0$).\\
And thus, if $\RR$ is nowhere hyperfinite then the partition of $U$ trivializes the decomposition.
\end{theorem}

Given the coincidence \cite{Gab02} of the first $\ell^2$-Betti number $\beta_1(\Gamma)$ of any countable group $\Gamma$ with the first $L^2$-Betti number $\beta_1(\RR_{\alpha},\mu)$ of the orbit equivalence relation defined by any free p.m.p. action $\Gamma\action{\alpha} (X,\mu)$ on the standard Borel space, 
and since non-amenability implies nowhere hyperfinite, we immediately get:
\begin{corollary}\label{cor: beta-1=0 -> Gamma is MFI}
Every non-amenable countable group $\Gamma$ with vanishing first $\ell^2$-Betti number $\beta_1(\Gamma)=0$ is measurably freely indecomposable.
\end{corollary}

\begin{question}
Produce examples {\MFI} groups with $\beta_1>0$.
\end{question}

\begin{question}
Characterize all the {\MFI} groups.
\end{question}

\bigskip
Let's say that a p.m.p. countable standard equivalence relation is \dfn{accessible} if it admits a free product decomposition $\RR=\ast_{j\in J} \RR_j$ into freely indecomposable subrelations. 
\begin{question}
Find/characterize  p.m.p. countable standard equivalence relation, with finite $\beta_1$, that are non-accessible.
\end{question}

\section{Bass-Serre Rigidity}

Suppose that $\Theta : G_1 \freeprod G_2 \longrightarrow G'_1 \freeprod G'_2$ is an isomorphism of groups where $G_i$ and $G'_j$ are freely indecomposable groups different from $\mathbb{Z}$. Since $\Theta(G_1)$ is a subgroup of $G'_1 \freeprod G'_2$, Kurosh's theorem implies that $\Theta(G_1)$ is a subgroup of a conjugate of $G'_1$ or $G'_2$. Up to a permutation of the indices, we assume that $\Theta(G_1)$ is a subgroup of a conjugate $\text{conj}(G'_1)$ of $G'_1$. Another use of Kurosh's theorem implies that $\Theta^{-1}(\text{conj}(G'_1))$ is a subgroup of a conjugate of $G_1$ or $G_2$. But since $\Theta^{-1}(\text{conj}(G'_1))$ contains $G_1$, we deduce equality with $G_1$, \textit{i.e.} $\Theta(G_1)=\text{conj}(G'_1)$ and in the same way that $\Theta(G_2)=\text{conj}(G'_2)$. This observation is the starting point of our main theorem:
\begin{theorem}\label{th: rigidity BS}
Let $\gauche{\RR}=\underset{{\gi}\in {\gI}}{*}\gauche{\RR}_{{\gi}}*\gauche{\TT}$
and $\droite{\RR}=\underset{{\di}\in {\dI}}{*}\droite{\RR}_{{\di}}*\droite{\TT}$ be p.m.p. standard equivalence relations decomposed into free products, where each factor $\gauche{\RR}_{{\gi}}$ and $\droite{\RR}_{{\di}}$ is \dfn{freely indecomposable} and aperiodic on its domain; and where $\gauche{\TT}$ and $\droite{\TT}$ are treeable.
If $\gauche{\RR}$ and $\droite{\RR}$ are SOE, via an isomorphism $\Theta : \gauche{V}\subset \locus{\gauche{\RR}}
\to \droite{V}\subset \locus{\droite{\RR}}$ then $\gauche{\RR}$ and $\droite{\RR}$ admit free product decompositions induced by countable slicings of the factors:
\begin{eqnarray}
	\forall \gi\in \gI, \hskip 5pt \locus{\gauche{\RR}_{\gi}}=\coprod\limits_{\gauche{\kk}\in {\gauche{\KK}(\gi)}} \gauche{X}_{\gauche{\kk}} &	\hfill \textrm{\ \  and\ \ \ }	\hfill &
\forall \di\in \dI, \hskip 5pt \locus{\droite{\RR}_{\di}}=\coprod\limits_{\droite{\kk}\in {\droite{\KK}(\di)}} \droite{X}_{\droite{\kk}} \\
	\gauche{\RR}=\underset{{\gi}\in {\gI}}{*}
	(
	\underset{\gauche{\kk}\in {\gauche{\KK}(\gi)}}{*}
	\gauche{\RR}_{\gi}\vert \gauche{X}_{\gauche{\kk}}	)	
	*\gauche{\TT} 
	&\hfill \textrm{\ \  and\ \ \ }	\hfill&
	\droite{\RR}=\underset{{\di}\in {\dI}}{*}
	(\underset{\droite{\kk}\in {\droite{\KK}(\di)}}{*}\droite{\RR}_{\di}\vert \droite{X}_{\droite{\kk}})
	*\droite{\TT}
\end{eqnarray}
for which there is a bijection $\theta: \coprod_{{\gi}\in {\gI}}\gauche{\KK}(\gi)\to \coprod_{{\di}\in {\dI}}\droite{\KK}(\di)$ between the index sets such that, denoting $\gauche{\SS}_{\gauche{\kk}}:=\gauche{\RR}_{\gi}\vert \gauche{X}_{\gauche{\kk}}$ and $\droite{\SS}_{\droite{\kk}}:=\droite{\RR}_{\di}\vert \droite{X}_{\droite{\kk}}$,
for each $\gauche{\kk}\in \coprod_{{\gi}\in {\gI}}\gauche{\KK}(\gi)$, the slices
$\gauche{\SS}_{\gauche{\kk}}$ and $\droite{\SS}_{\theta(\gauche{\kk})}$
 are SOE  via an isomorphism between subsets of the domains $\locus{\gauche{\SS}_{\gauche{\kk}}}=\gauche{X}_{\gauche{\kk}}$ and $\locus{\droite{\SS}_{\theta(\gauche{\kk})}}=\droite{X}_{\theta(\gauche{\kk})}$ of the shape $\droite{f} \Theta \gauche{f}$, where $\gauche{f}\in [[\gauche{\RR}]]$ and $\droite{f}\in [[\droite{\RR}]]$.
\end{theorem}
\begin{remark}
Recall that in case $\gauche{\RR}_{\gi}$ is ergodic, then it admits no non-trivial slicing. If all the $\gauche{\RR}_{\gi}$ and $\droite{\RR}_{\di}$ are ergodic, then Theorem~\ref{th: rigidity BS} establishes a bijection between the $\gauche{\RR}_{\gi}$ and the $\droite{\RR}_{\di}$
($\theta$ becomes a bijection between the sets of indices $\gI$ and $\dI$).
\end{remark}
The proof of the theorem will be given in Section~\ref{sect: proof th: rigidity BS}.
\section{Proofs}

\subsection{Preliminaries}\label{subsect: preliminaries}
We list some \dfn{properties} of $L^2$-Betti numbers of equivalence relations (see \cite{Gab02}) on the non atomic probability measure space $(X,\mu)$. Recall that the $L^2$-Betti numbers are defined with respect to an invariant probability measure \cite{Gab02}. In case a finite measure is invariant, one usually normalize it. Thus $\beta_q(\RR)$ stands for the $q$-th $\ell^2$-Betti number of $\RR$ on $\locus{\RR}$ with respect to the normalized probability measure $\frac{\mu\vert \locus{\RR}}{\mu(\locus{\RR})}$.
If $\locus{\RR}\subsetneq X$, the notation $\beta_q(\RR,\mu)$ means that we extend $\RR$ trivially outside $\locus{\RR}$ to compute $L^2$-Betti numbers according to the probability measure $\mu$. 
\begin{proposition}\label{prop: properties L2 Betti numbers}
The following holds:
\begin{enumerate}
\item $\beta_0(\RR)=\int_X \frac{d\mu(x)}{\#\RR(x)}$, the mean value of the inverse of the cardinal of the class of $x$, with the convention 
$\frac{1}{\infty}=0$. It follows that $\beta_0(\RR)\in [0,1]$.

\item The relation $\RR$ is trivial if and only if $\beta_0(\RR)=1$.\label{prop :beta0=0}

\item $\RR$ is aperiodic on $\locus{\RR}$ 
if and only if $\beta_0(\RR)=0$.\label{prop: beta-0=0 for infinite classes}

\item If $\beta_{1}(\RR)=0$ and $V\subset X$ satisfies $\mu(V)>0$, then $\beta_{1}(\RR{\vert V})=0$.\label{prop: beta-1=0 for restrictions}

\item For a free product: $\beta_{1}(\RR_1*\RR_2, \mu)-\beta_{0}(\RR_1*\RR_2,\mu)=\beta_1(\RR_1,\mu)-\beta_0(\RR_1,\mu)+\beta_1(\RR_2,\mu)-\beta_0(\RR_2,\mu)+1$.\label{prop: beta-1 of a free product}
\end{enumerate}
\end{proposition}
Proof: We use the notation of \cite{Gab02}.
If the space $(X,\mu)$ admits an $\RR$-invariant partition into non-negligible subsets $X=\coprod X_r$, then 
\begin{equation}
\beta_p(\RR,\mu)=\sum_{r} \mu(X_r) \beta_p(\RR\vert X_r, \frac{\mu\vert X_r}{\mu(X_r)})
\end{equation}
The Hilbert module $H$ one has to consider to define the $L^2$-Betti numbers may be decomposed into a direct sum $H=\oplus_r H_r$ according to the decomposition of $X$, and the normalization of the trace leads to $\dim_{\RR}H_r=\mu(X_r) \dim_{\RR\vert X_r} H_r$, and the formula. Up to partitioning the space into the Borel subsets $X_r$ where the classes have constant cardinal $r$, one may compute $\beta_0(\RR)$ under the assumption that the classes all have cardinal $r$. If $r=\infty$, then $\beta_0(\RR)=0$ (\cite[Prop. 3.15]{Gab02}). If $r$ is finite, there is a contractible $\RR$-complex containing only one point in each fiber. The computation is then immediate, $\beta_0(\RR)=\frac{1}{r}$. In general, 
\begin{equation}
\beta_0(\RR,\mu)=\sum_{r} \mu(X_r) \frac{1}{r}=\int_X \frac{d\mu(x)}{\#\RR(x)}\end{equation}
Properties \ref{prop :beta0=0}, \ref{prop: beta-0=0 for infinite classes} follow. Property~\ref{prop: beta-1=0 for restrictions} also follows, with the following standard lemma:
\begin{lemma}\label{lem: restr. of infinite classes=infinite}
If the classes of the standard finite measure preserving equivalence relation $\RR$ are almost all infinite, and $V\subset X$ is non-negligible, then almost all the classes of the restriction $\RR{\vert V}$ are infinite.
\end{lemma}
As for Property~\ref{prop: beta-1 of a free product} (which is a consequence of the Mayer-Vietoris sequence in the classical situations), we start by decomposing the space $X$ according to whether the classes of $\RR_1*\RR_2$ are finite or not. 
If the classes are almost all finite, then we are in a treeable situation and the formula is immediate by \cite[Cor. 3.23]{Gab02} and \cite[Th. IV.15]{Gab00a}. Assume now that the classes are infinite, then $\beta_0(\RR_1*\RR_2, \mu)=0$. 

Recall that the first two $L^2$-Betti numbers ($\beta_0$ and $\beta_1$) are computed by considering any $2$-dimensional simplicial $\RR$-complex $\Sigma$ with simply connected fibers, and any exhausting increasing family of uniformly locally bounded $\RR$-invariant sub-complexes $(\Sigma_t)_{t\in \Nmath}$. Then, one calculate the limits of the von Neumann dimensions (with respect to the von Neumann algebra associated with $\RR$ and the trace associated with $\mu$) for $*=0,1$:
\begin{equation}
\beta_{*}(\RR)=\lim\limits_{s\to \infty}\lim\limits_{t\to \infty, t\geq s} \underbrace{\dim_{\RR} \Image(H_{*}^{(2)}(\Sigma_s)\overset{J_{s,t}}{\to} H_{*}^{(2)}(\Sigma_t))}_{:=\nabla_{*}(\Sigma_s, \Sigma_t)}\label{eq: def of beta-* as lim lim}
\end{equation}
where ${J_{s,t}}$ is induced in homology by the inclusion $\Sigma_s\subset \Sigma_t$.

Let $\RR=\RR_1*\RR_2$.
Consider two such 2-dimensional complexes  and their exhaustions $\Sigma_1$, $(\Sigma_{1,t})_{t\in \Nmath}$ for $\RR_1$ and $\Sigma_2$, $(\Sigma_{2,t})_{t\in \Nmath}$ for $\RR_2$, and the induced $\RR$-complexes  $\SST_1$, $(\SST_{1,t})_{t\in \Nmath}$ and $\SST_2$, $(\SST_{2,t})_{t\in \Nmath}$, obtained by suspension \cite[sect. 5.2]{Gab02}. One may assume that $\SST_1^{0}=\SST_{1,t}^{0}=\SST_2^{0}=\SST_{2,t}^{0}\simeq \RR$. Then $\SST=\SST_1\cup\SST_2$ is a simply connected $2$-dimensional simplicial $\RR$-complex $\Sigma$, with associated exhausting sequence $\SST_t=\SST_{1,t}\cup\SST_{2,t}$. One may have arranged $\Sigma_1$ and $\Sigma_2$ so that $\beta_0(\SST_{1,t},\mu)=\beta_0(\RR_1,\mu)$, $\beta_0(\SST_{2,t},\mu)=\beta_0(\RR_2,\mu)$ and $\beta_0(\SST_t)=0$ for all $t\in \Nmath$ (\textit{i.e.} the connected components of $\SST_t$ are infinite).
Then 
\begin{eqnarray*}
\nabla_{1}({\widetilde{\Psi}}_s, {\widetilde{\Psi}}_t)-\nabla_{0}({\widetilde{\Psi}}_s, {\widetilde{\Psi}}_t)
&=& \dim_\RR C_1^{(2)}({\widetilde{\Psi}}_s)-1-\dim_\RR \Bigl( \Image \partial_2 C_2^{(2)}({\widetilde{\Psi}}_t)\cap C_1^{(2)}({\widetilde{\Psi}}_s)\Bigr)
\end{eqnarray*}
is valid for ${\widetilde{\Psi}}=\SST$, $\SST_1$ or $\SST_2$. Since the terms on the right hand side may be split 
\begin{eqnarray*}
C_1^{(2)}(\SST_s)&=&C_1^{(2)}(\SST_{1,s})\oplus C_1^{(2)}(\SST_{2,s})\\
\Image \partial_2 C_2^{(2)}(\SST_t)\cap C_1^{(2)}(\SST_s)&=&\Image \partial_2 C_2^{(2)}(\SST_{1,t})\cap C_1^{(2)}(\SST_{1,t})\\
&& \hskip 90pt \oplus \Image \partial_2 C_2^{(2)}(\SST_{2,t})\cap C_1^{(2)}(\SST_{2,t})
\end{eqnarray*}
it follows that 
\begin{multline}
\nabla_{1}(\SST_s, \SST_t)-\nabla_{0}(\SST_s, \SST_t)=\nabla_{1}(\SST_{1,s}, \SST_{1,t})-\nabla_{0}(\SST_{1,s}, \SST_{1,t})\\
+\nabla_{1}(\SST_{2,s}, \SST_{2,t})-\nabla_{0}(\SST_{2,s}, \SST_{2,t})+1
\end{multline}
and taking the limits, like in (\ref{eq: def of beta-* as lim lim}), leads to the formula of Property~\ref{prop: beta-1 of a free product}.\endofproof

\subsection{Proof of Theorem~\ref{th: trivial decomposition for beta-1=0}}
\label{sect: proof th: trivial decomposition for beta-1=0}
Consider $\RR$ with $\beta_{1}(\RR)=0$ and all the classes infinite.
We start assuming that $\RR$ decomposes as a free product of two factors $\RR=\RR_{1}*\RR_{2}$. 
Let $U_{i}$ be the union of the infinite $\RR_{i}$-classes, for $i=1,2$. 

\noindent
(1) We show that $\RR{\vert U_1}=\RR_{1}{\vert U_1}$, $\RR{\vert U_2}=\RR_{2}{\vert U_2}$ and 
 both $\RR_{2}{\vert U_1}$ and $\RR_{1}{\vert U_2}$ are trivial:
If $\mu(U_1)>0$, then the restrictions of $\RR$ and $\RR_1$ to $U_1$ satisfy respectively  $\beta_1(\RR{\vert U_1})=0$ (by \ref{prop: properties L2 Betti numbers} Property~\ref{prop: beta-1=0 for restrictions}) and $\beta_0(\RR_{1}{\vert U_1})=0$ (by \ref{prop: properties L2 Betti numbers} Property~\ref{prop: beta-0=0 for infinite classes}).
By Theorem~\ref{th: a la Kurosh}, $\RR{\vert U_1}=\RR_{1}{\vert U_1}*\SS$, where we have isolated the particular subrelation $\VV_{i_1}=\RR_{1}\cap\RR{\vert U_1}=\RR_{1}{\vert U_1}$  (Theorem~\ref{th: a la Kurosh}(\ref{item: R cap S-i is part of Kurosh decomposition})) and we have put all the other terms of the free product decomposition together to form $\SS$, which itself contains the other subrelation $\VV_{i_2}=\RR_{2}\cap\RR{\vert U_1}=\RR_{2}{\vert U_1}$. \\
By \ref{prop: properties L2 Betti numbers} Property~\ref{prop: beta-1 of a free product}, 
$\underbrace{\beta_1(\RR{\vert U_1})}_{=0}=\underbrace{\beta_{1}(\RR_{1}{\vert U_1})}_{\geq 0}+\underbrace{\beta_1(\SS)}_{\geq 0}+1-(\underbrace{\beta_0(\RR_{1}{\vert U_1})}_{=0}+\beta_0(\SS))$ so that $\beta_0(\SS)=1$, \textit{i.e.} (by Proposition~\ref{prop: properties L2 Betti numbers}, item~\ref{prop :beta0=0}) $\SS$ is trivial.
It follows that $\RR_{2}{\vert U_1}$ is trivial and the decomposition reduces to $\RR{\vert U_1}=\RR_{1}{\vert U_1}$. Symmetrically, if $\mu(U_2)>0$, $\RR_{1}{\vert U_2}$ is trivial and $\RR{\vert U_2}=\RR_{2}{\vert U_2}$.

\noindent
(2) 
We claim that $\mu(U_1\cap U_2)=0$, for otherwise, $(\RR_{1}{\vert U_1}){\vert U_1\cap U_2}$ the iterated restriction 
would have infinite classes (by Lemma~\ref{lem: restr. of infinite classes=infinite}). But it is also the trivial subrelation $(\RR_{1}{\vert U_2}){\vert U_1\cap U_2}$. 

\noindent
(3)  We claim that the partition $U=U_1\coprod U_2$ is $\RR{\vert U}$-invariant. 
If one of $U_1$, $U_2$ is a null set, 
this follows from part (1) of the proof, so that we may assume both are non-null. The partition is already $\RR_{1}{\vert U}$- and $\RR_{2}{\vert U}$-invariant.
As above, Theorem~\ref{th: a la Kurosh} gives a decomposition $\RR{\vert U}= \RR_{1}{\vert U}*\SS$, where $\SS$ contains $\RR_{2}{\vert U}$. The above parts (1) and (2) apply to this decomposition, with $U'_2$ the union of the infinite $\SS$-classes, in place of $U_2$, leading to $\mu(U_1\cap U_2')=0$. Observe that $U'_2$ contains the $\SS$-saturation of $U_2=U\setminus U_1$, so that $U'_2=U_2$ (a.s.) and $U_2$ is $\SS$-invariant. Being also $\RR_{1}{\vert U}$-invariant, $U_2$ ends up $\RR{\vert U}$-invariant. Symmetrically, $U_1$ is $\RR{\vert U}$-invariant. 

The four first points of Theorem~\ref{th: trivial decomposition for beta-1=0} have been proved for two factors.

\noindent
(4) If now $\RR=\RR_{1}*\RR_{2}*\cdots*\RR_{i}*\cdots$, we apply the above result after one factor $\RR_i$ has been isolated and the other ones have been glued together in an $\overline{\RR_{i}}$ leading to a decomposition $\RR=\RR_i*\overline{\RR_{i}}$. The union $\overline{U_i}$ of the infinite orbits of $\overline{\RR_{i}}$ contains all the $U_j$, for $j\not=i$.
We immediately deduce the four first points of Theorem~\ref{th: trivial decomposition for beta-1=0} in general. For instance, $U_i$ being $\RR{\vert U_i\cup \overline{U_i}}$-invariant is also invariant for the even more restricted equivalence relation $\RR{\vert (U_1\cup U_2\cup\cdots\cup U_i\cup\cdots)}$.

\noindent
(5) We conclude by proving that the restriction of $\RR$ to $W=X\setminus \RR.U$, the complement of the saturation of $U$, is hyperfinite as soon as $\mu(W)>0$. Observe that $W$ is $\RR$-invariant and that $\RR{\vert W}=\RR_{1}{\vert W}*\RR_{2}{\vert W}*\cdots*\RR_{i}{\vert W}*\cdots$. By definition of $U$, the restrictions $\RR_{i}{\vert W}$ are finite subrelations, and thus treeable, so that $\RR{\vert W}$ follows treeable and aperiodic. Since moreover $\beta_1(\RR{\vert W})=0$ (Proposition~\ref{prop: properties L2 Betti numbers} Property~\ref{prop: beta-1=0 for restrictions}) it is hyperfinite by 
Proposition~6.10 of \cite{Gab02}. If $\RR$ is nowhere hyperfinite, then $\mu(W)=0$ and $U$ is a complete section for $\RR$.
This completes the proof of Theorem~\ref{th: trivial decomposition for beta-1=0}.
\endofproof

\subsection{Proof of Theorem~\ref{th: rigidity BS}}
\label{sect: proof th: rigidity BS}
By Theorem~\ref{th: a la Kurosh restrict} and Proposition~\ref{prop: restr of FI is FI}, the proof reduces to the case where $\Theta:\locus{\gauche{\RR}} \to \locus{\droite{\RR}}$ is in fact an OE between
$\gauche{\RR}=\freeprod\limits_{{\gi}\in {\gI}}\gauche{\RR}_{{\gi}}*\gauche{\TT}$
and $\droite{\RR}=\freeprod\limits_{\di\in {\dI}}\droite{\RR}_{\di}*\droite{\TT}$.

a) Fix one ${\gi}\in {\gI}$ and define the subrelation 
\begin{equation}
\EE_{{\gi}}:=\Theta\gauche{\RR}_{{\gi}}\Theta^{-1}
\end{equation} 
of $\droite{\RR}$, the image of $\gauche{\RR}_{{\gi}}$ under $\Theta$. It admits a decomposition according to Theorem~\ref{th: a la Kurosh}:
\begin{equation}\label{eq: decomp. of E-p}
\EE_{{\gi}}=\freeprod\limits_{{\di}\in {\dI}}(\freeprod\limits_{\droite{\kk}\in \droite{\KK}(\gi,\di)} \droite{\VV}_{\droite{\kk}})*\droite{\TT}_{\gi}
\end{equation}
where for each $\droite{\kk}\in \droite{\KK}(\gi,\di)$, 
\begin{equation}
	\droite{\VV}_{\droite{\kk}}=\EE_{{\gi}}\cap {\droite{\psi}}^{-1}_{\droite	{\kk}}\ \droite{\RR}_{\di} \ \droite{\psi}_{\droite{\kk}}
\end{equation}
with $\droite{\psi}_{\kk}\in [[\droite{\RR}]]$ a partial isomorphism defined on $\locus{\droite{\VV}_{\droite{\kk}}}$. In particular, $\droite{\VV}_{\droite{\kk}}$ is inner conjugate in $\droite{\RR}$ with a subrelation of $\droite{\RR}_{\di}$. As for $\droite{\TT}_{\gi}$, it is a treeable subrelation, containing the treeable part given by Theorem~\ref{th: a la Kurosh} and the conjugates of subrelations of $\droite{\TT}$ (see Remark~\ref{rem: group the treeable part}).

Since $\EE_{{\gi}}$ is {\FI} just like $\gauche{\RR}_{{\gi}}$, this decomposition (\ref{eq: decomp. of E-p}) admits a trivializing partition
(Definition~\ref{def: trivialization of decompositon}).
Observe that the treeable part $\droite{\TT}_{{\gi}}$ cannot survive as a slice in the trivializing partition: its restriction to some $U_i$ would coincide with $\EE_{\gi}\vert U_i$, would be treeable \cite[Prop. II.6]{Gab00a} and (\FI) (Prop. \ref{prop: restr of FI is FI}); and thus smooth (Proposition~\ref{rem: treeable + FI -> smooth}), which is ruled out by the aperiodicity assumption on $\RR_{\gi}$.
The trivializing partition thus takes the form
\begin{equation}
	\coprod\limits_{{\di}\in {\dI}} \ \ \coprod\limits_{\droite{\kk}\in \droite{\KK}(\gi,\di)}	\UUU_{\droite{\kk}} 
\end{equation}
and induces the slicing (see Remark~\ref{Rem: sur la def de decomp. triviale} item~\ref{item: slicing induced by trivialization}. equation~(\ref{eq: slicing induced by trivialization})) affiliated with the partition of $\locus{\EE_{{\gi}}}$ into the $\EE_{{\gi}}$-saturations $\VVV_{\droite{\kk}} $ of $\UUU_{\droite{\kk}}$
\begin{eqnarray}\label{eq: slicing of Theta-1 R-i Theta}
\EE_{{\gi}}&=&\freeprod\limits_{{\di}\in {\dI}} \ \ {\freeprod\limits_{\droite{\kk}\in \droite{\KK}(\gi,\di)}} \EE_{{\gi}}\vert \VVV_{\droite{\kk}}\\
\label{eq: slicing of R-1}
	\gauche{\RR}_{\gi}=\Theta^{-1}\ \EE_{{\gi}}\ \Theta &=&\freeprod\limits_{{\di}\in {\dI}} \ \ {\freeprod\limits_{\droite	{\kk}\in \droite{\KK}(\gi,\di)}} 
	\gauche{\RR}_{{\gi}}	{\vert \Theta^{-1}\VVV_{\droite{\kk}}}
\end{eqnarray}
Moreover, by definition of the trivialization, for each $\droite{\kk}\in \droite{\KK}(\gi,\di)$ the restriction of $\EE_{{\gi}}$ to $\UUU_{\droite{\kk}}$ 
 satisfies:
\begin{eqnarray*}
\EE_{{\gi}}{\vert  \UUU_{\droite{\kk}}}
&=&{\droite{\VV}_{\droite{\kk}}{\vert  \UUU_{\droite{\kk}}}}\\
&=&\Bigl(\EE_{{\gi}}\cap \bigl({\droite{\psi}}^{-1}_{\droite	{\kk}}\ \droite{\RR}_{\di}
 \ \droite{\psi}_{\droite{\kk}}\bigr)\Bigr){\vert  \UUU_{\droite{\kk}}}\\
&=&{\EE_{{\gi}}{\vert  \UUU_{\droite{\kk}}}}\cap ({\droite{\psi}}^{-1}_{{\droite{\kk}}} \ \droite{\RR}_{\di} \  \droite{\psi}_{{\droite{\kk}}}){\vert  \UUU_{\droite{\kk}}}
\end{eqnarray*}
which means exactly
\begin{eqnarray}\label{eq: E-i subset of a unique factor}
\EE_{{\gi}}{\vert  \UUU_{\droite{\kk}}}&\subset& ({\droite{\psi}}^{-1}_{{\droite{\kk}}}\  \droite{\RR}_{\di} \  \droite{\psi}_{{\droite{\kk}}}){\vert  \UUU_{\droite{\kk}}} 
\end{eqnarray}

Since we did not yet use the properties of the $\droite{\RR}_{\di}$, let's raise what  we proved so far~:
\begin{proposition}
If $\SS=\underset{q\in Q}{*}{\SS}_{q}*{\TT}$  is p.m.p. and $\TT$ is treeable, and if 
${\mathcal{E}}$ is a subrelation that is \dfn{freely indecomposable} and aperiodic on its domain,
then there are (at most) countably many disjoint Borel subsets $U_{r}$ whose ${\mathcal{E}}$-saturations $V_{r}={\mathcal{E}} U_{r}$ form an ${\mathcal{E}}$-invariant partition $\locus{{\mathcal{E}}}=\coprod_{q\in Q} \coprod_{r\in R(q)} V_{r}$ with affiliated slicing 
${\mathcal{E}}=\underset{q\in Q}{*}\underset{r\in R(q)}{*} {{\mathcal{E}}}\vert V_{r}$ and such that  
for $r\in R(q)$, one has ${{\mathcal{E}}}\vert U_{r} \subset\psi_r ^{-1} \SS_{q} \psi_r\vert U_{r}$ for some $\psi_r\in [[\SS]]$.
\end{proposition}

b) Observe that the slicings (\ref{eq: slicing of Theta-1 R-i Theta}) of the factors
$\EE_{{\gi}}:=\Theta\gauche{\RR}_{{\gi}}\Theta^{-1}$ induce a corresponding slicing of $\droite{\RR}=\Theta\gauche{\RR}\Theta^{-1}$ and $\gauche{\RR}$:
\begin{eqnarray}
	\label{eq: free prod decomp of Theta R-1 Theta (-1)}
	\droite{\RR}&=&
	\freeprod_{{\gi}\in {\gI}}\bigl(\freeprod\limits_{{\di}\in {\dI}} \ \ {\freeprod\limits_{\droite	{\kk}\in \droite{\KK}(\gi,\di)}} \EE_{{\gi}}\vert \VVV_{\droite{\kk}}\bigr)\freeprod\Theta	\gauche{\TT}\Theta^{-1}
	\\
	\label{eq: refined free prod decomp of R-1}
	\gauche{\RR}&=&
	\freeprod_{{\gi}\in {\gI}}\bigl(\freeprod\limits_{{\di}\in {\dI}} \ \ {\freeprod\limits_{\droite	{\kk}\in \droite{\KK}(\gi,\di)}} 
	\gauche{\RR}_{{\gi}}	{\vert \Theta^{-1}\VVV_{\droite{\kk}}}\bigr)\freeprod\gauche{\TT}
\end{eqnarray}

c) The subrelation $({\droite{\psi}}^{-1}_{{\droite{\kk}}}\  \droite{\RR}_{\di} \  \droite{\psi}_{{\droite{\kk}}}){\vert  \UUU_{\droite{\kk}}}$
of $\droite{\RR}=\Theta\gauche{\RR}\Theta^{-1}$  
appearing in (\ref{eq: E-i subset of a unique factor}), for some $\gi\in \gI, \di\in \dI$ and $\droite{\kk}\in \droite{\KK}(\gi,\di)$ such that $\UUU_{\droite{\kk}}$ is non-negligible,  gets itself a free product decomposition with respect to (\ref{eq: free prod decomp of Theta R-1 Theta (-1)}), given by Theorem~\ref{th: a la Kurosh}: 
\begin{eqnarray}\label{eq: decomp R-2-k}
({\droite{\psi}}^{-1}_{{\droite{\kk}}}\  \droite{\RR}_{\di} \  \droite{\psi}_{{\droite{\kk}}}){\vert  \UUU_{\droite{\kk}}}
=\bigl(*_{l\in L} \WW_l  \bigr)*\TT_{\kk}
\end{eqnarray}
The point \ref{item: R cap S-i is part of Kurosh decomposition}. of Theorem~\ref{th: a la Kurosh} states that the particular term
\begin{eqnarray}
\WW_{l_0}&=&
({\droite{\psi}}^{-1}_{{\droite{\kk}}}\  \droite{\RR}_{\di} \  \droite{\psi}_{{\droite{\kk}}}){\vert  \UUU_{\droite{\kk}}}
\cap \EE_{{\gi}}{\vert  \VVV_{\kk}}
\end{eqnarray}
has to appear and from (\ref{eq: E-i subset of a unique factor}), we get:
\begin{eqnarray}
\WW_{l_0}&=&\EE_{{\gi}}{\vert  \UUU_{\droite{\kk}}}\label{eq: W-l=E-i}
\end{eqnarray}
On the other hand, $({\droite{\psi}}^{-1}_{{\droite{\kk}}}\  \droite{\RR}_{\di} \  \droite{\psi}_{{\droite{\kk}}}){\vert  \UUU_{\droite{\kk}}}$
is {\FI} since isomorphic with the restriction of the {\FI} relation $\droite{\RR}_{\di}$ to a non-null subset of its domain. As such, its decomposition (\ref{eq: decomp R-2-k}) admits a trivialization.

But the particular term $\EE_{{\gi}}{\vert  \UUU_{\droite{\kk}}}$ is nowhere smooth on its whole domain $\UUU_{\droite{\kk}}$, so that (Proposition~\ref{Prop: sur la def de decomp. triviale} item~\ref{item: R-i nowhere-smooth-> R=R-i}) 
 this term is the only one of the decomposition (\ref{eq: decomp R-2-k}); \textit{i.e.} (\ref{eq: E-i subset of a unique factor}) is an equality:
\begin{eqnarray}\label{eq: equality E-i= conj. R-k}
\EE_{{\gi}}{\vert  \UUU_{\droite{\kk}}}&=& ({\droite{\psi}}^{-1}_{{\droite{\kk}}}\  \droite{\RR}_{\di} \  \droite{\psi}_{{\droite{\kk}}}){\vert  \UUU_{\droite{\kk}}}\\
\Theta\ \  \gauche{\RR}_{{\gi}}{\vert \Theta^{-1}(\UUU_{\droite{\kk}})}\ \ \Theta^{-1} = (\Theta \RR_{{\gi}}\Theta^{-1}) {\vert  \UUU_{\droite{\kk}}}&=& 
{\droite{\psi}}^{-1}_{{\droite{\kk}}}\  \droite{\RR}_{\di}{\vert \droite{\psi}_{\droite{\kk}}(\UUU_{\droite{\kk}})} \  \droite{\psi}_{{\droite{\kk}}}
\end{eqnarray}

This shows that the map $\droite{\psi}_{\droite{\kk}}\Theta$ defines, for $\droite{\kk}\in \droite{\KK}(\gi,\di)$, an isomorphism 
\begin{equation}\label{eq: psi Theta gives an OE}
\droite{\psi}_{\droite{\kk}}\Theta: \gauche{\RR}_{{\gi}}{\vert \Theta^{-1}(\UUU_{\droite{\kk}})}\OE \droite{\RR}_{\di}{\vert \droite{\psi}_{\droite{\kk}}(\UUU_{\droite{\kk}})}
\end{equation}
and thus a SOE between the slicing term 
$\gauche{\RR}_{{\gi}}{\vert \Theta^{-1}(\VVV_{\droite{\kk}})}=\Theta^{-1} \ \EE_{{\gi}}{\vert  \VVV_{\droite{\kk}}}\ \Theta$ of (\ref{eq: slicing of R-1}) and $\droite{\RR}_{\di}{\vert \droite{\psi}_{\droite{\kk}}(\UUU_{\droite{\kk}})}$, the restriction of $\droite{\RR}_{\di}$ to $\droite{\psi}_{\droite{\kk}}(\UUU_{\droite{\kk}})$
\begin{equation}
\droite{\psi}_{\droite{\kk}}\Theta: \gauche{\RR}_{{\gi}}{\vert \Theta^{-1}(\VVV_{\droite{\kk}})}\SOE \droite{\RR}_{\di}{\vert \droite{\psi}_{\droite{\kk}}(\UUU_{\droite{\kk}})}
\end{equation}

\medskip

d) Fix a $\di\in\dI$. We will show that the family of sets 
\begin{equation}
\WWW(\droite{\kk}):=\droite{\psi}_{\droite{\kk}}(\UUU_{\droite{\kk}}) \hskip 30pt \textrm{\ for\ } \droite{\kk}\in \coprod_{\gi\in \gI}\droite{\KK}(\gi,\di)
\end{equation}
induces a slicing of $\droite{\RR}_{{\di}}$ (in particular the family is not empty), \textit{i.e.} we show that their $\droite{\RR}_{{\di}}$-saturation form a partition of $\locus{\droite{\RR}_{{\di}}}$.

d-1) We first show that their $\droite{\RR}_{{\di}}$-saturation intersect trivially.
Let $\droite{\kk}_1\in \droite{\KK}({\gi}_1,\di)$ and $\droite{\kk}_2\in \droite{\KK}({\gi}_2,\di) $ such that the $\droite{\RR}_{{\di}}$-saturations of $\WWW(\droite{\kk}_1)$ and $\WWW(\droite{\kk}_2)$ have a non-null intersection, \textit{i.e.} there is a partial isomorphism $\droite{\rho}\in [[\droite{\RR}_{{\di}}]]$ with (non-null) domain contained in $\WWW(\droite{\kk}_1)$ and target in $\WWW(\droite{\kk}_2)$. It follows that the partial isomorphism ${\droite{\psi}}^{-1}_{\droite{\kk}_2} \droite{\rho} \droite{\psi}_{\droite{\kk}_1}$ has non-null domain $A_1\subset \UUU_{\droite{\kk}_1}(\gi_1)$
and target $A_2\subset \UUU_{\droite{\kk}_2}(\gi_2)$ and conjugates 
$\EE_{{\gi}_1}{\vert  A_1}$ with $\EE_{{\gi}_2}{\vert  A_2}$.
But these subrelations are not smooth and appear as subrelations of factors of the free product decomposition (\ref{eq: free prod decomp of Theta R-1 Theta (-1)}). It follow from Lemma~\ref{lem: intersection de conj de facteurs} that they cannot belong to different factors, \textit{i.e.} $\droite{\kk}_1=\droite{\kk}_2$.

d-2) Consider now the invariant partition of $\locus{\droite{\RR}_{{\di}}}$ given by the $\droite{\RR}_{{\di}}$-saturations 
of the sets $\WWW(\droite{\kk})$, for $\droite{\kk}\in \coprod_{\gi\in \gI}\droite{\KK}(\gi,\di)$, and the complement $Z({\di})$ of their union in $\locus{\droite{\RR}_{{\di}}}$; and consider the affiliated slicing of $\droite{\RR}_{{\di}}$. We will show that the measure of $Z({\di})$ is zero.
We exchange the roles of $\gauche{\RR}$ and $\droite{\RR}$ after having further decomposed $\droite{\RR}$ thanks to the just above constructed slicing of $\droite{\RR}_{\di}$ affiliated with
$Z({\di})\coprod_{{\gi\in \gI}}\coprod_{\droite{\kk}\in\droite{\KK}(\gi,\di)}  \droite{\RR}_{{\di}}\WWW(\droite{\kk})$. We use $\Theta^{-1}$ and apply  the above steps a), b), c).
After a restriction of its domain, the slice $\droite{\RR}_{{\di}}{\vert Z({\di})}$ is conjugate with  one of the $\gauche{\RR}_{{\gi}}$ restricted to a Borel subset $Y$ of its domain, like in (\ref{eq: psi Theta gives an OE})
 via some $\gauche{\eta}_{1} \Theta^{-1}$, where $\eta_1\in [[\gauche{\RR}]]$. This restriction $\gauche{\RR}_{{\gi}}\vert Y$ has just been shown ((\ref{eq: psi Theta gives an OE}) again but in the direct sense) to be conjugate (up to an additional restriction) with a restriction of one of the $\droite{\RR}_{\droite{\qq}}{\vert \WWW(\droite{\qq})}$, for some $\droite{\qq}$ via some $\droite{\eta}_{2} \Theta$, with $\droite{\eta}_{2}\in [[\droite{\RR}]]$.
Since the composition $\droite{\eta}_{2} \Theta\gauche{\eta}_{1} \Theta^{-1}\in [[\droite{\RR}]]$, 
it follows that up to restricting to a non-negligible Borel subset, $\droite{\RR}_{{\di}}{\vert Z({\di})}$ is inner conjugate with a restriction of $\droite{\RR}_{\droite{\qq}}{\vert \WWW(\droite{\qq})}$, one of the factors in the decomposition of $\droite{\RR}$, which is different from $\droite{\RR}_{{\di}}{\vert Z({\di})}$ by definition of $Z({\di})$. The Lemma~\ref{lem: intersection de conj de facteurs} would imply that $\droite{\RR}_{{\di}}{\vert Z({\di})}$ is somewhere smooth, contrarily to the assumption that the orbits of the ${\droite{\RR}_{{\di}}}$ are all infinite on its domain. It follows that the measure of $Z({\di})$ is zero.

The families $\WWW(\droite{\kk})$ induce slicings of the factors $\droite{\RR}_{{\di}}$ leading to a refined free product decomposition of $\droite{\RR}$ whose (non treeable) terms are indexed by 
$\gauche{\KK}:=\coprod_{\gi\in \gI}\coprod_{\di\in \dI} \droite{\KK}(\gi,\di)$ and in a bijective SOE correspondence with those of the refined decomposition (\ref{eq: refined free prod decomp of R-1}) of $\gauche{\RR}$, via $\Theta$ and inner partial isomorphisms.
The slicing of $\gauche{\RR}_{\gi}$ we were after in Theorem~\ref{th: rigidity BS} is affiliated with the 
$A_{\droite{\kk}}:={\Theta^{-1} \VVV_{\droite{\kk}}}$ (the $\gauche{\RR}_{\gi}$-saturation of the ${\Theta^{-1} \UUU_{\droite{\kk}}}$),
for $\droite{\kk}\in \gauche{\KK}(\gi):=\coprod_{\di\in \dI} \droite{\KK}(\gi,\di)$, while the slicing of $\droite{\RR}_{\di}$ is affiliated with the $\droite{\RR}_{\di}$-saturation $B_{\droite{\kk}}:= \droite{\RR}_{\droite{\kk}}\WWW(\droite{\kk})=\droite{\RR}_{\droite{\kk}}\droite{\psi}_{\droite{\kk}}(\UUU_{\droite{\kk}})$, for $\droite{\kk}\in \droite{\KK}(\di):=\coprod_{\gi\in \gI} \droite{\KK}(\gi,\di)$, and the SOE is given by (\ref{eq: psi Theta gives an OE}). \endofproof

\subsection{Proof of the Corollaries}
\label{sect: proof of the Corollaries}
-- Proof of {Theorem~\ref{th: ME BS rigidity-intro}  of the introduction.}
By measure equivalence with free groups, we have treeable free p.m.p. actions $\Gamma_0\action{\alpha}X$ and $\Lambda_0\action{\beta}Y$. Considering a coupling $(\Omega,\nu)$ witnessing the measure equivalence of  Equation~(\ref{eq: ME between free prod}), the corresponding diagonal-action coupling on $\Omega\times X\times Y$, obtained by extending $\alpha,\beta$ trivially on the other factors, deliver SOE actions of $\freeprod_{i\in I}\Gamma_i\freeprod \Gamma_0$ and $\freeprod_{j\in J}\Lambda_j\freeprod \Lambda_0$ whose restrictions to 
$\Gamma_0$ and $\Lambda_0$ are treeable.
Theorem~\ref{th: ME BS rigidity-intro} then follows immediately from Theorem~\ref{th: rigidity BS}.\endofproof

\bigskip
\noindent
-- Proof of Corollary~\ref{Cor: ergodicity on one side}.
The ergodicity assumption on the $\gauche{\alpha}$ side prevents from any slicing for $\RR_{\gauche{\alpha}}$ in Theorem~\ref{th: rigidity BS} which gives nevertheless a bijection 
 $\theta: {\gI}\to \coprod_{{\di}\in {\dI}}\droite{\KK}(\di)$ ($\gauche{n}= \droite{n}$ follows) for which $\Theta$ induces a SOE between the terms $\RR_{\gauche{\alpha}\vert \gauche{\Gamma}_{\gi}}$ and $\RR_{\droite{\alpha}\vert \droite{\Gamma} _{\theta(\gi)}}$. The latter follows ergodic.
Under the moreover assumption, the free products $\freeprod_{\gi\in \gI}\gauche{\Gamma}_{\gi}$ and $\freeprod_{\di\in \dI}\droite{\Gamma}_{\di}$ have the same first $\ell^2$-Betti number, $\not=0, \infty$. Thus any SOE between them has to be an OE (\cite{Gab02}): $\Theta$ induces an OE between the ergodic subrelations $\RR_{\gauche{\alpha}\vert \gauche{\Gamma}_{\gi}}$ and $\RR_{\droite{\alpha}\vert \droite{\Gamma}_{\theta(\gi)}}$.\endofproof

\bigskip

\noindent
-- Corollary~\ref{cor: beta-p>0 --> OE} of the introduction is a specialization of the following.
\begin{theorem}\label{th: IME trivial -> OE}
Let $\Theta$ be a SOE between two actions $\gauche{\alpha}$ and $\droite{\alpha}$ as in Framework~\ref{hyp: general framework}. Assume that $\gauche{\Gamma}_1\ME \droite{\Gamma}_1$ with generalized index $1$, and that 
$\gauche{\Gamma}_1\notME \gauche{\Gamma}_{\gi},\droite{\Gamma}_{\di}$ 
for all $\gi,\di\not=1$.
Assume moreover that $I_{\mathrm{ME}}(\gauche{\Gamma}_1)=\{1\}$.
Then $\Theta$ is in fact an orbit equivalence and the restrictions to  $\gauche{\Gamma}_1$ and $\droite{\Gamma}_1$ are OE.
In particular, they have the same measure space of ergodic components, in particular the same families of measures of ergodic components (possibly with repetition).
\end{theorem}
Proof: Theorem~\ref{th: rigidity BS} applied to the SOE $\Theta$ between 
 $\RR_{\gauche{\alpha}}=\freeprod_{\gi\in \gI} \RR_{\gauche{\alpha}\vert \gauche{\Gamma}_\gi}\freeprod \RR_{\gauche{\alpha}\vert \gauche{\Gamma}_0}$ and 
 $\RR_{\droite{\alpha}}
 =\freeprod_{\di\in \dI}$
 $\RR_\droite{\alpha}\vert \droite{\Gamma}_{\di}
 \freeprod \RR_{\droite{\alpha}\vert \droite{\Gamma}_0}$ produces 
 slicings of $\RR_{\gauche{\alpha}\vert \gauche{\Gamma}_1}$ and $\RR_{\droite{\alpha}\vert \droite{\Gamma}_1}$ whose components are 
pairwise associated by $\theta$ and SOE via partial isomorphisms of the shape $\droite{f} \Theta \gauche{f}$ with $\droite{f}$, $\gauche{f}$ preserving respectively the measures $\gauche{\mu}, \droite{\mu}$. They all scale the measure by the same factor and may be assembled together in order to produce a global SOE between $\RR_{\gauche{\alpha}\vert \gauche{\Gamma}_1}$ and $\RR_{\droite{\alpha}\vert \droite{\Gamma}_1}$, with the same compression constant. The point being that all together the slices meet almost all their classes. Now $I_{\mathrm{ME}}(\gauche{\Gamma}_1)=\{1\}$  (for instance if some $\beta_p(\gauche{\Gamma}_1)\not=0, \infty$) and $\gauche{\Gamma}_1\MEa{1} \droite{\Gamma}_1$ imply that any SOE between free p.m.p. $\gauche{\Gamma}_1$- and $\droite{\Gamma}_1$-actions is in fact an OE. The compression constant equals $1$ and $\Theta$ is in fact an OE.
\endofproof

\nobreak
\noindent \textsc{Aur\'elien Alvarez \& Damien Gaboriau: \\
  Universit\'e de Lyon -- CNRS,\\
    ENS-Lyon, UMPA UMR 5669,\\
69364 Lyon
cedex 7,\\ FRANCE}

\noindent \texttt{aurelien.alvarez@umpa.ens-lyon.fr\\ damien.gaboriau@umpa.ens-lyon.fr}

\end{document}